\documentclass[11pt]{article}
\usepackage{amsfonts}
\usepackage{amssymb}
\usepackage{latexsym}
\usepackage{amscd}
\usepackage[centertags]{amsmath}
\usepackage{xypic}
\newtheorem{defn}{Definition}[section]
\newtheorem{thm}{Theorem}%[section]
\newtheorem{cor}[thm]{Corollary}
\newtheorem{lem}[thm]{Lemma}

\newtheorem{rem}{Remark}

\newcommand{\calN}{\mathcal{N}}

\newcommand{\rmd}{\mathrm{d}}

\newcommand{\Enum}{\mathbb{E}}
\newcommand{\Pnum}{\mathbb{P}}

\newcommand{\Rnum}{\mathbb{R}}

\newcommand{\Znum}{\mathbb{Z}}
\newcommand{\Nnum}{\mathbb{N}}
\newcommand{\Var}{\mathrm{Var}\,}

\newcommand{\lam}{\lambda}

\newcommand{\qed}{\hfill $\Box$}
\newcommand{\wh}{\widehat}
\newcommand{\wt}{\widetilde}

\begin{document}
%===========================================================================================
\title{Range-Renewal Processes: SLLN, Power Law and Beyonds}

\author{Xin-Xing Chen$^{1}$, Jian-Sheng Xie$^{2}$\footnote{Corresponding author. E-mail: jiansheng.xie@gmail.com}
and Jiangang Ying$^{2}$\\
{\footnotesize 1. Department of Mathematics, Shanghai Jiao Tong University, Shanghai 200240, China}\\
{\footnotesize 2. School of Mathematical Sciences, Fudan University, Shanghai 200433, China}}

\date{}
\maketitle
%===========================================================================================
\begin{abstract}
Given $n$ samples of a regular discrete distribution $\pi$, we prove in this article first a
serial of SLLNs results (of Dvoretzky and Erd\"{o}s' type) which implies a typical power law
when $\pi$ is heavy-tailed. Constructing a (random) graph from the ordered $n$ samples, we can
establish other laws for the degree-distribution of the graph. The phenomena of small world 
is also discussed.
%===========================================================================================
\end{abstract}
%===========================================================================================
%%===========================================================================================
%\begin{keyword}
%%===========================================================================================
%%% keywords here, in the form: keyword \sep keyword
%
%%% MSC codes here, in the form: \MSC code \sep code
%%% or \MSC[2008] code \sep code (2000 is the default)
%Range-renewal \sep SLLN \sep Regular Varying \sep Power Law \sep Small World
%%===========================================================================================
%\end{keyword}
%%===========================================================================================

% \linenumbers

\section{Introduction}\label{sec:1}
%===========================================================================================================
%Suppose that a spider selects and moves to a vertex according to a
%given probability law $\pi$ at each step, independent of its
%previous choice; the trace it left behind forms a random (network)
%graph. We would like to know the number of vertices of the graph
%formed in the first $n$ steps. Such a problem can be modeled
%mathematically in the following way:
Let $\xi:=\{\xi_n: n \geq 1\}$ be a symbol sequence (with certain distribution law $\pi$)
and let $R_n$ be the number of distinct elements among the first $n$ elements of the
process $\xi$; Let's call $\{R_n: n \geq 1\}$ the range-renewal process of $\xi$. We
would like to investigate the growth rate of $R_n$ (and other related quantities such as
$R_{n, \,k}$ to be defined later) in $n$. Of course, for non-triviality, we should assume
that the process $\xi$ is in fact infinitely (but discretely) valued. Our main concern in
this article is the independent and identically distributed case, meaning that $\xi :=
\{\xi_n: n \geq 1\}$ is i.i.d.. Therefore the model can be interpreted as the following:
Suppose that a spider selects and moves to a vertex according to a given probability law
$\pi$ at each step, independent of its previous choice; the trace it left behind forms a
random graph (or network). In this setting, we prove first the strong law of large numbers
(of Dvoretzky and Erd\"{o}s' type) $\lim \limits_{n \to +\infty} R_n /\Enum R_n=1$; then
under some mild and natural assumptions on the common distribution $\pi$, we prove a
sequence of SLLNs for other related quantities which implies a power law and other laws. As
we know, there are thousands of works (see e.g. \cite{BA99, AH00, B09, AB02} and references
therein and thereafter) concerning power laws and other related things of all kinds of
random graph models; A significant part of our current work was influenced and inspired by
them.

Let's explain how we come to study such a problem and also to find
the approach which we would present later.

In the autumn of 2010 the second author reported the classic result
of Dvoretzky and Erd\"{o}s \cite{DE50} in a seminar at Fudan
University and was fascinated by their neat and beautiful result
\begin{equation}\label{eq: 1.1}
\lim_{n \to +\infty} \frac{R_n}{\Enum R_n}=1
\end{equation}
for simple symmetric random walks (abbr. \textbf{SSRW}) on $\Znum^d$ with $d
\geq 2$ of course; here $R_n$ denotes the number of cites visited by
the random walk in the first $n$ steps. Let's call $R_n$ the
\textbf{range-renewnal} at time $n$. Erd\"{o}s and Taylor had further discussions
about \textbf{SSRW} on $\Znum^d$ \cite{ET}, which closely relates to the work \cite{DE50}.

We then try to find out the more recent results concerning
$R_n$ for maybe more general processes. Chosid-Isaac and Athreya
\cite{CI78} \cite{CI80} \cite{Athreya85} obtained the limit
\begin{equation}\label{eq: 1.1'}
\lim_{n \to +\infty} \frac{R_n}{n}
\end{equation}
being zero for irreducible positively recurrent Markov chain or null
recurrent Markov chain under a suitable integrable condition.
Derriennic \cite{Derriennic} extended Dvoretzky-Erd\"{o}s' result to
random walk (based on stationary distributions) on discrete groups.
He showed that the limit in (\ref{eq: 1.1'}) always exists almost
surely. Furthermore if the random walk on the group is recurrent,
the limit is zero; otherwise the limit is just the escape rate as
Dvoretzky-Erd\"{o}s' result says for SSRW on $\Znum^d$ with $d \geq
3$. The central limit theorem for $R_n$ (for SSRW on $\Znum^d$) can
be found in Jain and Pruitt \cite{JP71} \cite{JP74} ($d \geq 3$) and
Le Gall \cite{LeGall86} ($d =2$). Law of the iterated logarithm for
$R_n$ (for SSRW on $\Znum^d$) are discussed by Bass and Kumagai
\cite{BK02} ($d=3$) and Jain and Pruitt \cite{JP72} ($d \geq 4$).
More discussions on null recurrent or transient Markov chains can be
found in \cite{Glynn85} \cite{GS92} \cite{Hamana97} \cite{Hamana98}
\cite{SV05} \cite{Vallois96} \cite{VT97} and references therein.

As it is already seen, there are fruitful results concerning $R_n$
for null recurrent or transient Markov chains. By contrast, in
general there are relatively few results concerning $R_n$ for
positive recurrent Markov chains (or more general, stationary
processes) with infinite denumerable states. All we now know in general is that
$R_n \uparrow \infty$ and (see e.g. \cite{CI78})
\begin{equation}\label{eq: 1.3}
\lim_{n \to +\infty} \frac{R_n}{n}=0.
\end{equation}
But what would be the accurate order of $R_n$ tending to $+\infty$?
This is an interesting and important problem. This problem has not
been investigated even for i.i.d sequence (the simplest Markov chain
model and the basic assumption in statistics and sampling) ever
since the publication of Dvoretzky and Erd\"{o}s' result for random
walk \cite{DE50} in 1950. This is our original motivation of this
research.

And the research of the above problem leads us to the current simple
and interesting criteria (which could be traced implicitly back to Dvoretzky and
Erd\"{o}s \cite{DE50}): Let $\displaystyle S_n :=\sum_{k=1}^n
\eta_k$ be a sum of non-negative random variables. Suppose $\Enum
S_n \to +\infty$ and $\sup\{\Enum \eta_n: n \geq 1\}<+\infty$.
Assume furthermore that we have the following estimation
\begin{equation}
\label{eq: variation-estimate'} \mathrm{Var} \, (S_n) \leq C \cdot
(\Enum S_n)^{2-\delta}
\end{equation}
for some positive $C, \delta$ and all $n$, or even more weakly
\begin{equation}
\label{eq: variation-estimate''} \mathrm{Var} \, (S_n) \leq C \cdot
(\Enum S_n)^2/(\log \Enum S_n)^{1+\delta},
\end{equation}
then we can derive the following \textbf{strong law of large
numbers} (abbr. \textbf{SLLN})
\begin{equation}
\label{eq: SLLN} \lim_{n \to \infty} \frac{ S_n}{\Enum S_n}=1.
\end{equation}
Making use of the above approach, we give an almost complete answer to the question
proposed for i.i.d. models in the proceeding paragraph. We note here that, one year
later after finishing the proof of our main results for i.i.d. models in 2011, the
second author (along with J. Wu) discovered a similar structure in continued fractions
\cite{WX12}. The criteria just mentioned plays a crucial role in their proof as well as
in this paper; It also indicates the significance of our current research.

%The criteria just mentioned plays a crucial role in their proof. It
%also indicates the significance of our current research.
%This argument works for simple symmetric random walk on $\Znum^d$
%with $d \geq 3$ (and even $d=2$ with the corresponding estimation in
%\cite{DE50}; the proof there implies this criteria though they did
%not state this explicitly); see Example \ref{examp: 1} in Appendix
%B. It's interesting that, by this approach, we can give a
%new and simple proof of Borel's SLLN (see Example \ref{examp: 3} in
%Appendix B). The most important is, by this argument, we have proved a
%sequence of SLLNs (and power laws and other results under mild
%conditions) for i.i.d symbol sequence, which is a surprisingly new
%and beautiful result. With our approach it is possible to build
%similar SLLNs for more general positive recurrent Markov chains (see
%our Example \ref{examp: 4} in Appendix B) other than i.i.d case or
%for null recurrent/transient Markov chains other than R.W. on
%$\Znum^d$ (with $d=2$ or $d \geq 3$), provided that one can obtain
%estimation (\ref{eq: variation-estimate'}) for the model under
%investigation. We note here that, one year later after finishing the
%proof of our main results for i.i.d. models in 2011, the second author (along
%with J. Wu) discovered a similar structure in continued fractions \cite{WX12}.
%The criteria just mentioned plays a crucial role in their proof. It
%also indicates the significance of our current research.

The paper is organized in the following way. Section \ref{sec:2}
devotes to the presentation of the main settings, assumptions,
main results and related discussions, where small world
phenomena is discussed. In Section \ref{sec:3}
we present some necessary estimations for our model. Section \ref{sec:4}
is devoted to the proof of the main Theorems \ref{thm: main-non-critical},
\ref{thm: main-sup-critical} and \ref{thm: main-sub-critical}. %We
%discuss in Appendix A the relation between the range-renewal speed and
%the entropy of the distribution. In Appendix B, we discuss how our SLLN
%approach can be applied in some concrete examples, where an interesting
%new proof of Borel's SLLN is presented based on our SLLN approach.
%===========================================================================================
\section{Main Settings, Assumptions, Main Results and Related Discussions}\label{sec:2}
%===========================================================================================
%===========================================================================================
\subsection{Main Settings and Assumptions}
%===========================================================================================
Let $\{\xi_n: n \geq 1\}$ be a sequence of i.i.d. random variables
with common distribution $\pi$ which, for non-triviality and simplicity, is assumed to be
supported on the natural numbers set $\Nnum$. We denote by $R_n$ the number of distinct
values of $\xi_k, k=1, \cdots, n$, i.e.
\begin{equation}\label{eq:renewalrank-0}
R_n :=\#\{\xi_k: 1 \leq k \leq n\}.
\end{equation}
More importantly, given the random sequence $\xi_1, \cdots, \xi_n$,
we can obtain a finite directed random graph $G_n :=(V_n, E_n)$,
where $V_n$ is the set of vertices $\xi_1, \cdots, \xi_n$ and $E_n$
is the set of directed edges $\xi_i \to \xi_{i+1}$ (with starting
vertex $\xi_i$ and ending vertex $\xi_{i+1}$) for $i=1, \cdots,
n-1$; the induced undirected graph would be denoted by
$\widehat{G}_n :=(V_n, \widehat{E}_n)$ where $\widehat{E}_n$ is the
induced edges set. Clearly, $R_n$ is just the size of $V_n$, i.e.,
$R_n=\#(V_n)$. We call $\{R_n: n \geq 1\}$ the \textbf{range-renewal
process} with respect to the original process $\{\xi_n: n \geq 0\}$.

For latter use, put
\begin{equation}
N_n (x) :=\sum_{k=1}^n 1_{\{\xi_k=x\}},
\end{equation}
which is the number of visiting times (visiting intensity) at vertex
$x$ of the random sequence up to time $n$. Then put for each $\ell
\geq 1$
\begin{equation}
R_{n, \, \ell} :=\sum_x 1_{\{N_n (x)=\ell\}}, \quad R_{n, \, \ell+} :=\sum_x 1_{\{N_n (x) \geq \ell\}};
\end{equation}
Thus $R_{n, \, \ell}$ (resp. $R_{n, \, \ell+}$) is the number of distinct states which have been visited at
exactly (resp. at least) $\ell$ times in the first $n$ steps. Obviously
$$
R_n=R_{n,\, 1+}=\sum_{\ell=1}^n R_{n, \, \ell}, \quad R_{n, \,
\ell+}=\sum_{k=\ell}^n R_{n, \, k}.
$$
We will call the above numbers $R_{n, \, \ell}, R_{n, \, \ell+}$
(and other related numbers) \textbf{visiting intensity statistics}.

Define for any $x,y$
\begin{equation}
d_n (x, y) :=\sum_{i=1}^{n-1} 1_{\{\xi_i=x, \, \xi_{i+1}=y\}}.
\end{equation}
This is the visiting intensity of the edge $x \to y$ in the graph
$G_n$. Define also
\begin{equation}
D_n (x) :=\sum_y 1_{\{d_n (x,y) \geq 1\}}.
\end{equation}
This is the out-degree of vertex $x$ in graph $G_n$. We then define
for any $1 \leq k \leq \ell$
\begin{equation}
\wt{R}_{n, \, k, \, \ell} := \sum_x 1_{\{N_{n-1} (x)=\ell, \, D_n
(x)=k\}},\quad \wt{R}_{n, \, k} := \sum_x 1_{\{D_n (x)=k\}}.
\end{equation}
Clearly $\wt{R}_{n, \, k, \, \ell}$ is the number of vertices which
have out-degree$=k$ in graph $G_n$ but visiting intensity$=\ell$ in
graph $G_{n-1}$; $\wt{R}_{n, \, k}$ is the number of vertices with
out-degree$=k$ in graph $G_n$.

In order to investigate the undirected graphs $\wh{G}_n$, define
similarly
$$
\wh{d}_n (x, y) :=\sum_{i=1}^{n-1} 1_{\{\xi_i=x, \, \xi_{i+1}=y
\hbox{ or } \xi_i=y, \, \xi_{i+1}=x\}}.
$$
This is the visiting intensity of the undirected edge
$\overline{xy}$ in the graph $\wh{G}_n$. Define also
\begin{equation}
\wh{D}_n (x) :=\sum_y 1_{\{\wh{d}_n (x,y) \geq 1\}}.
\end{equation}
This is the degree of vertex $x$ in graph $\wh{G}_n$. Put for any
interval $\Delta \subset \Rnum^+$
$$
N_x (\Delta):=\sum_{k \in \Delta} 1_{\{\xi_k=x\}}.
$$
We then define for any $1 \leq k \leq 2 \ell$
\begin{equation}
\wh{R}_{n, \, k, \, \ell} := \sum_x 1_{\{N_x ([2, n-1])=\ell, \,
\wh{D}_n (x)=k\}}, \quad \wh{R}_{n, \, k} := \sum_x 1_{\{\wh{D}_n (x)=k\}}.
\end{equation}
Clearly $\wh{R}_{n, \, k, \, \ell}$ is the number of vertices which
have degree$=k$ in graph $\wh{G}_n$ but visiting intensity$=\ell$ in
graph $\wh{G}_{[2, n-1]}$, the graph induced by the sequence
$\xi_2, \cdots, \xi_{n-1}$; $\wh{R}_{n, \, k}$ is the number of
vertices with out-degree$=k$ in graph $\wh{G}_n$. We will call the
above numbers $\wt{R}_{n, \, k, \, \ell}, \wt{R}_{n, \, k}$ (resp.
$\wh{R}_{n, \, k, \, \ell}, \wh{R}_{n, \, k}$) (and other related
quantities) \textbf{out-degree statistics} (resp. \textbf{degree
statistics}).

For simplicity of the following discussion, we will assume that
\begin{itemize}
  \item[(\textbf{C0})] the distribution $\pi$ satisfies $\pi_{_1} \geq
\pi_{_2} \geq \cdots$ and $\pi_n>\pi_{n+1}>0$ for all large enough $n$.
We will denote
\begin{equation}\label{eq:def-phi}
\varphi (x)=1/\pi_x
\end{equation}
(where $x \in \Nnum$). Moreover, we would assume that the function
$\varphi (x)$ (see eq. (\ref{eq:def-phi})) is in fact continuously
defined for all $x \in [1, \infty)$ so that $\varphi (x)$ is
strictly increasing in $x$ for large enough $x$; this means that the
inverse function $\varphi^{-1} (x)$ exists for large enough $x$.
\end{itemize}

Furthermore we would require that the inverse function $\varphi^{-1}
(x)$ satisfies the first (if $\gamma \in (0,1)$) or the first and
second (if $\gamma=0$) or the first and third (if $\gamma=1$)
assumptions listed below (with suitable dominations in the related
limits):
\begin{itemize}
  \item[(\textbf{C1})] For some $\gamma \in [0,1]$, we have
\begin{equation}\label{eq: C2'-1}
\lim_{n \to \infty} \frac{\varphi^{-1} (n \lam)}{\varphi^{-1}
(n)}=\lam^\gamma, \quad \forall \lam \in (0, \infty);
\end{equation}
  \item[(\textbf{C1$^{\prime}$})] For $\gamma=0$, the function $\varphi_0 (x) :=\varphi^{-1}
(e^x)$ is $C^1$-smooth such that its derivative $\varphi_0^\prime
(x)$ satisfies: for any fixed $b \in \Rnum$
\begin{equation}\label{eq: C2'-1.0}
\lim_{x \to +\infty} \frac{\varphi_0^\prime (x+b)}{\varphi_0^\prime
(x)}=1.
\end{equation}
  \item[(\textbf{C1$^{\prime\prime}$})] For $\gamma=1$, there exist an increasing function
$\psi: [0, \infty) \to \Rnum^+$ and a continuous $L^1$-function $g:
[0, \infty) \to \Rnum^+$  (its $L^1$-norm will be denoted by
$\|g\|_1$) such that
\begin{equation}\label{eq: C2'-1.1}
\psi (n) \uparrow \infty \; \hbox{ and } \lim_{n \to \infty}
\frac{\varphi^{-1} (n \cdot e^{\lam \psi (\log n)})}{\varphi^{-1}
(n) \cdot e^{\lam \psi (\log n)}}=g (\lam) , \quad \forall \lam>0.
\end{equation}
\end{itemize}

The assumptions (\textbf{C1$^{\prime}$}) and
(\textbf{C1$^{\prime\prime}$}) are proposed to give a further
treatment for the critical cases of $\gamma=0$ and $\gamma=1$
respectively. %The first critical case of $\gamma=0$ will be called
%\textbf{sub-critical} with the second critical case of $\gamma=1$
%named \textbf{sup-critical}.
In general, for all assumptions, we would require reasonable
dominations on the limits so that Lebesgue's \textbf{Dominated
Convergence Theorem} can be applied in our discussion; but for
simplicity of the presentation, these requirements are {\it not
stated explicitly} in the above assumptions.

%===========================================================================================
\begin{defn}
%===========================================================================================
A distribution $\pi$ on $\Nnum$ is called \textbf{non-critical}, if
it satisfies assumption (\textbf{C1}) with $0<\gamma<1$ and suitable
domination in the related limit (\ref{eq: C2'-1}). It is called
\textbf{sub-critical}, if it satisfies assumptions
(\textbf{C1})+(\textbf{C1$^{\prime}$}) with $\gamma=0$ and suitable
dominations in the related limits (\ref{eq: C2'-1}) and (\ref{eq:
C2'-1.0}). It is called \textbf{sup-critical}, if it satisfies
assumptions (\textbf{C1})+(\textbf{C1$^{\prime\prime}$}) with
$\gamma=1$ and suitable dominations in the related limits (\ref{eq:
C2'-1}) and (\ref{eq: C2'-1.1}). Sometimes we would write the index
$\gamma=\gamma (\pi)$ to indicate its dependence on $\pi$.
%===========================================================================================
\end{defn}
%===========================================================================================
%===========================================================================================
\begin{defn}
%===========================================================================================
A distribution $\pi$ on $\Nnum$ is called \textbf{regular}, if $\pi$
is either non-critical, or sub-critical, or sup-critical.
%===========================================================================================
\end{defn}
%===========================================================================================

%===========================================================================================
\begin{rem}
%===========================================================================================
Eq. (\ref{eq: C2'-1}) is just the definition of regularly varying function of index $\gamma$,
which is originally introduced by Karamata \cite{Karamata}. See, e.g., \cite[pp. 321--324]{I-K}
or \cite[pp. 241--250]{Feller} for the definition and related properties. Eq. (\ref{eq: C2'-1.0})
is also related to regularly varying functions.
%===========================================================================================
\end{rem}
%===========================================================================================

%===========================================================================================
\begin{rem}
%===========================================================================================
It's easy to see that:
\begin{itemize}
  \item[{\rm (1)}] if $\pi_x=\frac{C}{x^\alpha} \cdot (1+o(1))$
with $1<\alpha<\infty, C>0$, then it satisfies (\textbf{C1}) with
$\gamma=1/\alpha \in (0,1)$;
  \item[{\rm (2)}] if $\pi_x=C \cdot e^{-a \cdot x} \cdot
(1+o(1))$ with $C, a>0$, then it satisfies (\textbf{C1}) with
$\gamma=0$ and (\textbf{C1$^{\prime}$}) with $\varphi_0 (x)=x^a
\cdot (1+o(1))/C$;
  \item[{\rm (3)}] if $\pi_x=\frac{C}{x \cdot (\log x)^\beta} \cdot
(1+o(1))$ with $\beta>1$ and $C>0$, then it satisfies (\textbf{C1})
with $\gamma=1$ and (\textbf{C1$^{\prime\prime}$}) with $\psi (x)=x$
and $g (\lam)=\frac{1}{(1+\lam)^\beta}$.
\end{itemize}
%===========================================================================================
\end{rem}
%===========================================================================================

%===========================================================================================
\subsection{Main results}
%===========================================================================================
Our main tool in this article is the following lemma.
%===========================================================================================
\begin{lem}\label{lem: SLLN-approach}
%===========================================================================================
Let $\displaystyle S_n :=\sum_{k=1}^n \eta_k$ be a sum of
non-negative random variables $\{\eta_n: n \geq 1\}$. Suppose $\Enum
S_n \to +\infty$ and $M:=\sup\{\Enum \eta_n: n \geq 1\}<+\infty$.
Furthermore we have the following estimation
\begin{equation}
\mathrm{Var} \, (S_n) \leq C \cdot (\Enum S_n)^{2-\delta}
\end{equation}
for some positive $C, \delta$ and all $n$; or even more weakly
\begin{equation}
\mathrm{Var} \, (S_n) \leq C \cdot (\Enum S_n)^2/(\log \Enum
S_n)^{1+\delta},
\end{equation}
then $\displaystyle \lim_{n \to \infty} \frac{ S_n}{\Enum S_n}=1$ almost surely.
%===========================================================================================
\end{lem}
%===========================================================================================
Since the proof of the above lemma is in fact contained implicitly in \cite{DE50} (with slight
mordifications) and is indeed an easy application of Borel-Cantelli lemma by noting that the
integer part of $\Enum S_n/M$ can run over all positive integers as $n$ running over all
positive integers, it is left to the readers as an exercise.

Our first main result is the following.
%===========================================================================================
\begin{thm}\label{thm: SLLN4Rn}
%===========================================================================================
For i.i.d. model with common distribution $\pi$ supported on a
countably infinite atoms, we always have $\displaystyle \lim_{n \to \infty} \frac{R_n}{\Enum R_n}
=1$ almost surely along with the following formulae
\begin{eqnarray}
\label{eq:expect} \Enum R_n &=& \sum_x [1- (1-\pi_x)^n] =: E (n),\\
\Var (R_n) &\leq& \Enum R_n.
\end{eqnarray}
%===========================================================================================
\end{thm}
%===========================================================================================

For regular distribution $\pi$, we have the following results.
%===========================================================================================
\begin{thm}\label{thm: main-non-critical}
%===========================================================================================
Assume $\pi$ to be \textbf{non-critical} (with $0<\gamma=\gamma
(\pi)<1$).
\begin{itemize}
  \item[{\rm (1)}] For any $\ell \geq 1$
\begin{equation}\label{eq:structure4noncritical-1}
\lim_{n \to \infty} \frac{R_{n, \, \ell}}{R_n} = r_\ell (\gamma)
:=\frac{\gamma \cdot \Gamma (\ell -\gamma)}{\ell! \cdot \Gamma
(1-\gamma)}=\frac{\gamma \cdot \prod_{j=1}^{\ell-1}
(j-\gamma)}{\ell!}
\end{equation}
holds true almost surely. It's clear that
\begin{equation}
r_\ell (\gamma)=\frac{\gamma}{\Gamma (1-\gamma)} \cdot
\ell^{-(1+\gamma)} \cdot [1+O (\frac{1}{\ell})]
\end{equation}
as $\ell \to \infty$, which is a \textbf{power law}. (\ref{eq:structure4noncritical-1})
is also equivalent to
\begin{equation}\label{eq:structure4noncritical-1.1}
\lim_{n \to \infty} \frac{R_{n, \, \ell}}{R_{n, \, \ell+}} =
\frac{\gamma}{\ell},
\end{equation}
meaning that the proportion of the relatively ``new" vertices at level
$\ell$ is approximately $\gamma/\ell$; this is a kind of average
escape rate (at level $\ell$); In the case of SSRW on $\Znum^d$ with $d \geq 3$, the
limit in (\ref{eq:structure4noncritical-1.1}) is always $\gamma_d$, the
usual escape rate, see e.g. \cite[p. 220]{Revesz}. Furthermore,
\begin{eqnarray}
\lim_{n \to \infty} \frac{R_{n, \, \ell}}{\Enum R_{n, \, \ell}} &=&
1, \\
\Enum R_n &=& \Gamma (1-\gamma) \cdot \varphi^{-1} (n) \cdot [1+o
(1)],\\
\Enum R_{n, \, \ell} &=& \frac{\gamma \cdot \Gamma (\ell
-\gamma)}{\ell!} \cdot \varphi^{-1} (n) \cdot [1+o (1)].
\end{eqnarray}
  \item[{\rm (2)}] Moreover, for each $1 \leq k \leq \ell$
\begin{eqnarray}
\label{eq:structure4noncritical-2'} \lim_{n \to \infty}
\frac{\wt{R}_{n, \, k, \, \ell}}{R_{n, \, \ell}} =S_{k, \, \ell}
(\pi) := \sum_{\#\{x_1, \cdots, x_\ell\}=k} \prod_{j=1}^\ell \pi_{x_j}, \\
\label{eq:structure4noncritical-2} \lim_{n \to \infty}
\frac{\wt{R}_{n, \, k}}{R_n} =f_k (\pi) := \sum_{\ell=k}^{\infty}
r_\ell (\gamma) \cdot S_{k, \, \ell} (\pi),
\end{eqnarray}
also hold true almost surely. And
\begin{eqnarray}
\lim_{n \to \infty} \frac{\wt{R}_{n, \, k, \, \ell}}{\Enum
\wt{R}_{n, \, k, \, \ell}} &=& 1, \quad \lim_{n \to \infty} \frac{\wt{R}_{n, \, k}}{\Enum \wt{R}_{n, \,
k}} =1, \\
\Enum \wt{R}_{n, \, k, \, \ell} &=& S_{k, \, \ell} (\pi) \cdot
\frac{\gamma \cdot \Gamma (\ell-\gamma)}{\ell!} \cdot
\varphi^{-1} (n) \cdot [1+o (1)],\\
\Enum \wt{R}_{n, \, k} &=& f_k (\pi) \cdot \Gamma (1-\gamma) \cdot
\varphi^{-1} (n) \cdot [1+o (1)].
\end{eqnarray}
  \item[{\rm (3)}] Similarly, for each $1 \leq k \leq 2 \ell$
\begin{eqnarray}
\label{eq:structure4noncritical-2'-undirected} \lim_{n \to \infty}
\frac{\wh{R}_{n, \, k, \, \ell}}{R_{n, \, \ell}} &=& S_{k, \, 2
\ell} (\pi), \\
\label{eq:structure4noncritical-2-undirected} \lim_{n \to \infty}
\frac{\wh{R}_{n, \, k}}{R_n} &=& \wh{f}_k (\pi) :=
\sum_{\ell=\lfloor (k+1)/2 \rfloor}^{\infty} r_\ell (\gamma) \cdot S_{k, \, 2
\ell} (\pi).
\end{eqnarray}
\end{itemize}
%===========================================================================================
\end{thm}
%===========================================================================================

%===========================================================================================
\begin{thm}\label{thm: main-sup-critical}
%===========================================================================================
Assume $\pi$ to be \textbf{sup-critical} (with $\gamma=\gamma
(\pi)=1$).
\begin{itemize}
  \item[{\rm (1)}] We have almost surely
\begin{eqnarray}
\lim_{n \to \infty} \frac{R_{n, \, 1}}{R_n} &=& 1, \qquad
\lim_{n \to \infty} \frac{R_{n, \, \ell}}{R_n} = 0, \quad \ell
\geq 2, \\
\lim_{n \to \infty} \frac{R_{n, \, \ell}}{R_{n, 2+}} &=& \frac{1}{\ell
\cdot (\ell-1)}, \quad \ell \geq 2;
\end{eqnarray}
The last limited ratio is still a \textbf{power law}; Equivalently
\begin{equation}
\lim_{n \to \infty} \frac{R_{n, \, \ell}}{R_{n, \ell+}} =
\frac{1}{\ell}, \quad \ell \geq 1.
\end{equation}
Furthermore,
\begin{eqnarray}
\lim_{n \to \infty} \frac{R_{n, \, \ell}}{\Enum R_{n, \, \ell}} &=&
1, \\
\Enum R_n &=& \|g\|_1 \cdot \varphi^{-1} (n) \cdot \psi (\log n) \cdot [1+o (1)],\\
\Enum R_{n, \, 1} &=& \|g\|_1 \cdot \varphi^{-1} (n) \cdot \psi
(\log n) \cdot [1+o (1)],\\
\Enum R_{n, \, \ell} &=& \frac{\varphi^{-1} (n)}{\ell \cdot
(\ell-1)} \cdot [1+o (1)], \quad \ell \geq 2,\\
\Enum R_{n, \, \ell+} &=& \frac{\varphi^{-1} (n)}{\ell-1} \cdot [1+o
(1)], \quad \ell \geq 2.
\end{eqnarray}
  \item[{\rm (2)}] Moreover, $\wt{R}_{n, \, 1, \, 1}=R_{n-1, \, 1}$
and
\begin{eqnarray}
\lim_{n \to \infty} \frac{\wt{R}_{n, \, 1}}{R_n} &=&
\lim_{n \to \infty} \frac{\wt{R}_{n, \, 1, \, 1}}{R_n}=1,\\
\lim_{n \to \infty} \frac{\wt{R}_{n, \, k}}{R_n} &=& 0, \quad k \geq
2
\end{eqnarray}
almost surely. Also a re-scaling yields for each $1 \leq k \leq
\ell$
\begin{eqnarray}
\label{eq:structure4supcritical-2'} \lim_{n \to \infty}
\frac{\wt{R}_{n, \, k, \, \ell}}{R_{n, \, \ell}} &=& S_{k, \, \ell}
(\pi) = \sum_{\#\{x_1, \cdots, x_\ell\}=k} \prod_{j=1}^\ell
\pi_{x_j}\\
\label{eq:structure4supcritical-2} \lim_{n \to \infty}
\frac{\wt{R}_{n, \, k}}{R_{n, \, 2+}} &=& f_k (\pi) :=
\sum_{\ell=k}^{\infty} \frac{S_{k, \, \ell} (\pi)}{\ell \cdot
(\ell-1)}, \quad k \geq 2\\
\lim_{n \to \infty} \frac{\wt{R}_{n, \, 1, \ell}}{R_{n, \, 2+}} &=&
\frac{S_{1, \, \ell} (\pi)}{\ell \cdot (\ell-1)}, \quad \ell \geq 2
\end{eqnarray}
almost surely. And
\begin{equation}
\lim_{n \to \infty} \frac{\wt{R}_{n, \, k, \, \ell}}{\Enum
\wt{R}_{n, \, k, \, \ell}} = 1, \quad
\lim_{n \to \infty} \frac{\wt{R}_{n, \, k}}{\Enum \wt{R}_{n, \, k}}
= 1
\end{equation}
almost surely. Moreover, for any $1 \leq k \leq \ell$ with $\ell
\geq 2$
\begin{eqnarray}
\Enum \wt{R}_{n, \, k, \, \ell} &=& \frac{S_{k, \, \ell} (\pi)}{\ell
\cdot (\ell-1)} \cdot \varphi^{-1} (n) \cdot [1+o (1)],\\
\Enum \wt{R}_{n, \, k} &=& f_k (\pi) \cdot \Gamma (1-\gamma) \cdot
\varphi^{-1} (n) \cdot [1+o (1)], \quad k \geq 2.
\end{eqnarray}

  \item[{\rm (3)}] Similarly, for each $1 \leq k \leq \ell$
\begin{eqnarray}
\lim_{n \to \infty} \frac{\wh{R}_{n, \, 1}}{R_n} &=& 1,
\lim_{n \to \infty} \frac{\wh{R}_{n, \, k}}{R_n} = 0, \quad k \geq 2\\
\lim_{n \to \infty} \frac{\wh{R}_{n, \, k, \, \ell}}{R_{n, \, \ell}}
&=& S_{k, \, 2 \ell} (\pi)\\
\lim_{n \to \infty} \frac{\wh{R}_{n, \, k}}{R_{n, \, 2+}} &=&
\wh{f}_k (\pi) := \sum_{\ell=\lfloor (k+1)/2 \rfloor}^{\infty} \frac{S_{k, \, 2
\ell} (\pi)}{\ell \cdot (\ell-1)}, \quad k \geq 2\\
\lim_{n \to \infty} \frac{\wh{R}_{n, \, 1, \ell}}{R_{n, \, 2+}} &=&
\frac{S_{1, \, 2 \ell} (\pi)}{\ell \cdot (\ell-1)}, \quad \ell \geq
2.
\end{eqnarray}
\end{itemize}
%===========================================================================================
\end{thm}
%===========================================================================================

%===========================================================================================
\begin{thm}\label{thm: main-sub-critical}
%===========================================================================================
Assume $\pi$ to be \textbf{sub-critical} (with $\gamma=\gamma
(\pi)=0$).
\begin{itemize}
  \item[{\rm (1)}] We have $\displaystyle \lim_{n \to \infty} \frac{R_{n, \, \ell}}{R_n} = 0, \quad \ell \geq
1$ almost surely. Furthermore,
\begin{eqnarray}
\Enum R_n &=& \varphi^{-1} (n) \cdot [1+o (1)],\\
\Enum R_{n, \, \ell} &=& \frac{\varphi_0^{\prime} (n)}{\ell} \cdot [1+o (1)], \quad \ell \geq 1,\\
\Enum R_{n, \, \ell+} &=& \varphi^{-1} (n) \cdot [1+o (1)], \quad
\ell \geq 2.
\end{eqnarray}
  \item[{\rm (2)}] Moreover, for each $k \geq 1$ $\displaystyle \lim_{n \to \infty} \frac{\wt{R}_{n, \, k}}{R_n} =0$
almost surely. And for each $1 \leq k \leq \ell$
\begin{eqnarray}
\Enum \wt{R}_{n, \, k, \, \ell} &=& \frac{S_{k, \, \ell}
(\pi)}{\ell}
\cdot \varphi_0^{\prime} (n) \cdot [1+o (1)],\\
\Enum \wt{R}_{n, \, k} &=& \wt{f}_k (\pi) \cdot \varphi_0^{\prime}
(n) \cdot [1+o (1)],\\
\Enum \wt{R}_{n, \, k+} &=& \varphi^{-1} (n) \cdot [1+o (1)],
\end{eqnarray}
where $\displaystyle \wt{f}_k (\pi):=\sum_{\ell=k}^\infty \frac{S_{k, \, \ell}
(\pi)}{\ell}$.
  \item[{\rm (3)}] Similarly, for each $k \geq 1$
$\displaystyle \lim_{n \to \infty} \frac{\wh{R}_{n, \, k}}{R_n} = 0$.
\end{itemize}
%===========================================================================================
\end{thm}
%===========================================================================================

When $\pi$ is not regular, we have the counter-examples.
%===========================================================================================
\begin{thm}\label{thm: main-counterexample}
%===========================================================================================
For any $0 \leq \gamma_1<\gamma_2 \leq 1$, there exists some
distribution $\pi$ which does not satisfy the assumption
(\textbf{C1}), such that there exists an increasing sequence
$\{n_j\}_{j=1}^{\infty} \subset \Nnum$ with
\begin{equation}\label{eq: counterexample} \lim_{j \to \infty} \frac{R_{n_{_{2j-1}}, 1}}{R_{n_{_{2j-1}}}}
=\gamma_1, \quad \lim_{j \to \infty} \frac{R_{n_{_{2j}},
1}}{R_{n_{_{2j}}}} =\gamma_2
\end{equation}
almost surely. %This means that, if assumption (\textbf{C1}) fails,
%the related limits in the above three theorems may also fail to
%exist almost surely.
%===========================================================================================
\end{thm}
%===========================================================================================
For the convenience of the reader, we present a proof of the above
theorem based on our theorem \ref{thm: main-non-critical} right now;
the proofs for theorems \ref{thm: main-non-critical}--\ref{thm:
main-sub-critical} would be given in the successive sections.

%===========================================================================================
\noindent \textit{Proof of Theorem \ref{thm:
main-counterexample}.\;}
%===========================================================================================
For simplicity, here we only construct a counter-example for
$0<\gamma_1<\gamma_2<1$; the other cases can be treated similarly.

For any distribution $\pi$ on $\Nnum$, we would denote by
$\Pnum_{\pi}$ the probability measure for the i.i.d. sequence of
$\{\xi_n: n \geq 1\}$ with common distribution $\pi$. For
simplicity, we will write $\displaystyle V_n:=\{\xi_i: i=1, \cdots, n\}$.

First we put $\alpha_1:=1/\gamma_1, \alpha_2 :=1/\gamma_2$ and
define a distribution $\pi^{(1)}$ on $\Nnum$ by
$$
\pi_x^{(1)} :=\frac{1}{Z_1 \cdot x^{\alpha_1}}, \quad x \in \Nnum,
$$
where $Z_1:=\sum_x \frac{1}{x^{\alpha_1}}$ is the normalizing
constant. In view of Theorem \ref{thm: main-non-critical},
$\Pnum_1:=\Pnum_{\pi^{(1)}}$-almost surely we clearly have
$\displaystyle \lim_{n \to \infty} \frac{R_{n, 1}}{R_n} =\gamma_1$.
Thus there exist large enough $n_1, m_1 \geq 1$ such that
$$
\Pnum_1 \Bigl( |\frac{R_{n_{_1}, 1}}{R_{n_{_1}}} -\gamma_1| \leq
\frac{\gamma_1}{2}, V_{n_{_1}} \subset [1, m_1) \Bigr) \geq
1-\frac{1}{2}
$$
and $\displaystyle \sum_{x \geq m_1} \pi^{(1)}_x \leq \frac{1}{2}$.
Then we adjust $\pi^{(1)}$ into a new distribution $\pi^{(2)}$ on
$\Nnum$ so that
$$
\pi^{(2)}_x :=\left\{
\begin{array}{rcl}
\pi^{(1)}_x, &\quad& 1 \leq x<m_1\\
\frac{1}{Z_2 \cdot x^{\alpha_2}}, &\quad& x \geq m_1
\end{array}
\right.
$$
In view of Theorem \ref{thm: main-non-critical},
$\Pnum_2:=\Pnum_{\pi^{(2)}}$-almost surely we have
$\displaystyle \lim_{n \to \infty} \frac{R_{n, 1}}{R_n} =\gamma_2$.
Thus there exist large enough $n_2 >n_1, m_2 >m_1$ such that
$$
\Pnum_2 \Bigl( |\frac{R_{n_{_2}, 1}}{R_{n_{_2}}} -\gamma_2| \leq
\frac{\gamma_2}{2 \cdot 2}, V_{n_{_2}} \subset [1, m_2) \Bigr) \geq
1-\frac{1}{2^2}
$$
and $\displaystyle \sum_{x \geq m_2} \pi^{(2)}_x \leq \frac{1}{2 \cdot 2}$.
Inductively, suppose we have already constructed a distribution
$\pi^{(2k)}$ with index $\gamma (\pi^{(2k)})=\gamma_2$, we are in a
position to construct a new distribution $\pi^{(2k+1)}$ with index
$\gamma (\pi^{(2k+1)})=\gamma_1$. Clearly, in view of Theorem
\ref{thm: main-non-critical}, $\Pnum_{2k} :=\Pnum_{\pi^{(2k)}}$-almost surely
we have $\displaystyle \lim_{n \to \infty} \frac{R_{n, 1}}{R_n} =\gamma_2$.
Thus there exist large enough $n_{2k} >n_{2k-1}, m_{2k} >m_{2k-1}$ such that,
if we write
\begin{eqnarray*}
A_k &=& \{|\frac{R_{n_{2k-1}, 1}}{R_{n_{2k-1}}} -\gamma_1| \leq
\frac{\gamma_1}{2 \cdot (2k-1)}, \; V_{n_{2k-1}} \subset [1,
m_{2k-1})\}, \\
B_k &=& \{|\frac{R_{n_{2k}, 1}}{R_{n_{2k}}} -\gamma_2| \leq
\frac{\gamma_2}{2 \cdot (2k)}, \; V_{n_{2k}} \subset [1, m_{2k})\},
\end{eqnarray*}
then $\displaystyle \Pnum_{2k} \Bigl(B_k\Bigr) \geq 1-\frac{1}{2^{2k}}$
and $\displaystyle \sum_{x \geq m_{2k}} \pi^{(2k)}_x \leq \frac{1}{2 \cdot (2k)}$.
Then we adjust $\pi^{(2k)}$ into $\pi^{(2k+1)}$ as the following:
$$
\pi^{(2k+1)}_x :=\left\{
\begin{array}{rcl}
\pi^{(2k)}_x, &\quad& 1 \leq x<m_{2k}\\
\frac{1}{Z_{2k+1} \cdot x^{\alpha_1}}, &\quad& x \geq
m_{2k}
\end{array}
\right.
$$
Also we can adjust $\pi^{(2k+1)}$ into $\pi^{(2k+2)}$ in the same
spirit. And finally we obtain a distribution $\pi^*$ on $\Nnum$:
$\displaystyle \pi^*_x :=\lim_{n \to \infty} \pi^{(n)}_x$.

It is easy to see that for each $k \geq 1$ $\pi^*_x=\pi^{(k)}_x, \quad \forall x<m_k$.
Put $\Pnum_*:=\Pnum_{\pi^*}$. It is clear that for each $k \geq 1$
\begin{eqnarray*}
\Pnum_* \Bigl( A_k \Bigr) &=& \Pnum_{2k-1} \Bigl( A_k \Bigr) \geq
1-\frac{1}{2^{2k-1}},\\
\Pnum_* \Bigl( B_k \Bigr) &=& \Pnum_{2k} \Bigl( B_k \Bigr) \geq
1-\frac{1}{2^{2k}}.
\end{eqnarray*}
Therefore we clearly have (\ref{eq: counterexample}) $\Pnum_*$-almost surely.
%===========================================================================================
\qed
%===========================================================================================

%===========================================================================================
\subsection{Discussion: Small World Phenomena}\label{sec:2.3}
%===========================================================================================
Now let's consider the diameter $L_n$ of the induced undirected
graph $\wh{G}_n$. Clearly, $\displaystyle L_n=\sup_{\overline{xy} \in E (\wh{G}_n)} L (x,y;
\wh{G}_n)$, where $L (x,y; \wh{G}_n)$ denotes the smallest length of a path
between vertices $x$ and $y$ in $\wh{G}_n$. Conditioned on
$\xi_1=x_0$ for some fixed $x_0$, denote by $T_1>1, T_2, \cdots$ the
successive times that the process $\{\xi_n: n \geq 1\}$ visits the
state $x_0$; put $\tau_k:=T_k-T_{k-1}$ with $T_0:=1$. It is obvious
that $\displaystyle L_n \leq \max\{\tau_k: 1 \leq k \leq N_n (x_0)+1\}$.

Noting that $\{\tau_k\}_{k=1}^\infty$ is an i.i.d. sequence with
common distribution $\Pnum (\tau_1=m)=\pi_{x_0} (1-\pi_{x_0})^{m-1},
m=1,2, \cdots$, we have for $\wh{\tau}_n:=\max\{\tau_k: 1 \leq k
\leq n\}$
$$
\Pnum (\wh{\tau}_n \leq \ell)=[1-(1-\pi_{x_0})^\ell]^n
$$
and
$$
\Enum \wh{\tau}_n=\frac{-\log n}{\log (1-\pi_{x_0})}+O (\log \log n),
\mathrm{Var\,} (\wh{\tau}_n) \leq O ((\log n) \log \log n),
$$
which implies $\displaystyle \lim_{n \to +\infty} \frac{\wh{\tau}_n}{\log n}
=\frac{-1}{\log (1-\pi_{x_0})}$
almost surely in view of Lemma \ref{lem: SLLN-approach}. Furthermore,
we would have
$$
\lim_{n \to +\infty} \frac{\max \limits_{1 \leq k \leq N_n (x_0)+1}
\tau_k}{\log n}=\frac{-1}{\log (1-\pi_{x_0})}
$$
since $N_n (x_0)/n \to \pi_{x_0}$ as $n \to +\infty$. Therefore we
always have
$$
\varlimsup_{n \to +\infty} \frac{L_n}{\log n} \leq C:=\inf_x [\frac{-1}{\log (1-\pi_{x})}]
=\frac{-1}{\log (1-\sup \limits_x \pi_{x})}<+\infty
$$
for any distribution $\pi$. Noting that, when $\pi$ is non-critical
or sup-critical, the size of the graph $\wh{G}_n$ is $R_n$ with
$\log R_n=O (\log n)$. Therefore, in these two case, we always have
$$
\varlimsup_{n \to +\infty} \frac{L_n}{\log R_n} \leq
C^\prime<+\infty
$$
for some constant $C^\prime$ which is a small world phenomena.

The accurate order of $L_n$ of the graph $\wh{G}_n$ seems to be much lower than
$\log R_n$. We guess still a rough bound $(\log R_n)^\gamma$ (where $\gamma=\gamma (\pi)$) such that
$$
\varlimsup_{n \to +\infty} \frac{L_n}{(\log R_n)^\gamma} <+\infty \hbox{ almost surely }
$$
for non-critical distribution $\pi$; but the calculation would be rather hard.
The accurate order of $L_n$ is surely an interesting (and even harder) open problem.

\section{Preliminary Estimates}\label{sec:3}
%===========================================================================================
%===========================================================================================
\subsection{Expectation Estimates for Visiting Intensity Statistics}
%===========================================================================================
Let's first start the estimate of $E (n)=\Enum R_n$ as a heat-up.
%===========================================================================================
\begin{lem}\label{lem: estimate-1}
%===========================================================================================
The function
\begin{equation}\label{eq: funct4expect}
E (z) :=\sum_x [1- (1-\pi_x)^z]
\end{equation}
is analytic on the complex plane with $\mathrm{Re\,} (z)>0$ (i.e.,
complex numbers with positive real parts). Furthermore, its $k$-th
derivative can be written as
$$
E^{(k)} (z) :=\frac{\mathrm{d}^k}{\mathrm{d} z^k} E (z)=-\sum_x
(1-\pi_x)^z \cdot [\log (1-\pi_x)]^k,
$$
which is still analytic on the complex plane with $\mathrm{Re\,}
(z)>0$.
\begin{itemize}
  \item[{\rm (i)}] For non-critical distribution $\pi$ (with $0<\gamma<1$),
\begin{eqnarray}
\label{eq: E-1} E (n) &=& \Gamma (1-\gamma) \cdot \varphi^{-1} (n) \cdot \Bigl[1+o
(1)\Bigr],\\
\label{eq: E-2} E^{(k)} (n) &=& \frac{(-1)^{k-1} \cdot \gamma \Gamma (k-\gamma)}{n^k}
\cdot \varphi^{-1} (n) \cdot \Bigl[1+o (1)\Bigr], \; k \geq 1;
\end{eqnarray}

  \item[{\rm (ii)}] For sub-critical distribution $\pi$ (with $\gamma=0$),
\begin{eqnarray}
\label{eq: E-1'} E (n) &=& \varphi^{-1} (n) \cdot \Bigl[ 1+o (1)\Bigr],\\
\label{eq: E-2'} E^{(k)} (n) &=& \frac{(-1)^{k-1} \cdot \Gamma (k)}{n^k} \cdot
\varphi_0^\prime (n) \cdot \Bigl[1+o (1)\Bigr], \; k \geq 1;
\end{eqnarray}
  \item[{\rm (iii)}] For sup-critical distribution $\pi$ (with $\gamma=1$),
\begin{eqnarray}
\label{eq: E-1"} E (n) &=& \|g\|_1 \cdot \varphi^{-1} (n) \cdot \psi (\log n) \cdot
\Bigl[1+o (1)\Bigr], \\
\label{eq: E-2"} E^{(1)} (n) &=& \frac{\|g\|_1}{n} \cdot \varphi^{-1} (n) \cdot \psi
(n) \cdot \Bigl[1+o (1)\Bigr],\\
\label{eq: E-3"} E^{(k)} (n) &=& \frac{(-1)^{k-1} \cdot (k-2)!}{n^k} \cdot
\varphi^{-1} (n) \cdot \Bigl[1+o (1)\Bigr], \quad k \geq 2.
\end{eqnarray}
\end{itemize}
%===========================================================================================
\end{lem}
%===========================================================================================
%===========================================================================================
\noindent \textit{Proof. \;}
%===========================================================================================Since
Since $\displaystyle R_n=\sum_x 1_{\{N_n (x) \geq 1\}}$, we clearly have
$$
\Enum R_n=\sum_x \Pnum(N_n (x) \geq 1)=\sum_x [1- (1-\pi_x)^n].
$$
That's why we study the function $E (z)$ defined by (\ref{eq: funct4expect}).

We only give a proof for {\rm (i)}; the other cases can be proved similarly. By our assumption,
$\varphi (x)$ is strictly increasing in $x$ (at least for large enough $x$). Thus the discrete
sum $\sum_x [1-(1-\pi_x)^n]$ can be approximated by the integral $\int_1^\infty [1-
(1-\frac{1}{\varphi (x)})^n] \rmd x$ with the error term bounded by 1. By our assumptions on
$\pi$, it is not hard to prove eq. (\ref{eq: E-1}). Eq. (\ref{eq: E-1}) can be proved similarly.
%===========================================================================================
\qed
%===========================================================================================
%===========================================================================================
\begin{cor}\label{cor: sup-critical-order}
%===========================================================================================
For the sup-critical case (with $\gamma=1$), we know that the
functions $\varphi$ and $\psi$ in (\textbf{C1$^{\prime\prime}$})
should satisfy
\begin{equation}\label{eq: sup-critical-order}
\lim_{n \to \infty} \frac{\log \varphi^{-1} (n)}{\log n}=1, \quad
\lim_{n \to \infty} \frac{\log \psi (\log n)}{\log n}=0.
\end{equation}
%===========================================================================================
\end{cor}
%===========================================================================================

Since we want to study $R_{n, \, \ell}$ for all $\ell \geq 1$, we
calculate out that
\begin{equation}
E_\ell (n) :=\Enum R_{n, \, \ell}=\sum_x C_n^\ell \cdot \pi_x^\ell
\cdot (1-\pi_x)^{n-\ell},
\end{equation}
where $C_n^\ell=\frac{n!}{\ell! \cdot (n-\ell)!}$. So let's write
\begin{equation}
S_\ell (n) :=\sum_x \pi_x^\ell \cdot (1-\pi_x)^{n-\ell}.
\end{equation}
Also we put
$$
E_{\ell+} (n) :=\Enum R_{n, \, \ell+}=\sum_{j=\ell}^n E_j (n).
$$
In the same spirit we have the following lemma.
%===========================================================================================
\begin{lem}\label{lem: estimate-2}
%===========================================================================================
For fixed integer $\ell \geq 1$, the function
$$
S_\ell (z) :=\sum_x \pi_x^\ell \cdot (1-\pi_x)^{z-\ell}
$$
is analytic on the right-half complex plane. Furthermore, its $k$-th
derivative can be written as
$$
S_\ell^{(k)} (z) :=\frac{\mathrm{d}^k}{\mathrm{d} z^k} S_\ell (z)
=\sum_x \pi_x^\ell (1-\pi_x)^{z-\ell} \cdot [\log (1-\pi_x)]^k,
$$
which is still analytic on the right-half complex plane.
\begin{itemize}
  \item[{\rm (i)}] For non-critical distribution $\pi$ (with $0<\gamma<1$),
\begin{eqnarray*}
S_\ell (n) &=&  \frac{\gamma \Gamma (\ell-\gamma)}{n^\ell}
\cdot \varphi^{-1} (n) \cdot \Bigl[1+o (1)\Bigr], \; \ell \geq 1,\\
S_\ell^{(k)} (n) &=& \frac{(-1)^{k} \cdot \gamma \Gamma
(k+\ell-\gamma)}{n^{k+\ell}} \cdot \varphi^{-1} (n) \cdot \Bigl[1+o
(1)\Bigr], \; k, \ell \geq 1
\end{eqnarray*}
and
\begin{eqnarray*}
E_\ell (n) &=& \frac{\gamma \Gamma
(\ell-\gamma)}{\ell!} \cdot \varphi^{-1} (n)\cdot \Bigl[1+o (1)\Bigr], \; \ell \geq 1, \\
E_{\ell+} (n) &=& \frac{\Gamma (\ell-\gamma)}{(\ell-1)!} \cdot
\varphi^{-1} (n) \cdot \Bigl[1+o (1)\Bigr], \ell \geq 2;
\end{eqnarray*}
  \item[{\rm (ii)}] For sub-critical distribution $\pi$ (with $\gamma=0$),
\begin{eqnarray*}
S_\ell (n) &=& \frac{(\ell-1)!}{n^\ell} \cdot \varphi_0^\prime (\log
n) \cdot \Bigl[1+o (1) \Bigr], \; \ell \geq 1,\\
S_\ell^{(k)} (n) &=& \frac{(-1)^{k} \cdot (k+\ell-1)!}{n^{k+\ell}}
\cdot \varphi_0^\prime (\log n) \cdot \Bigl[1+o (1)\Bigr], \; k,
\ell \geq 1
\end{eqnarray*}
and
\begin{eqnarray*}
E_\ell (n) &=& \frac{\varphi_0^\prime (\log n)}{\ell} \cdot \Bigl[1+o (1)\Bigr], \; \ell \geq 1, \\
E_{\ell+} (n) &=& \varphi^{-1} (n) \cdot \Bigl[1+o (1)\Bigr], \ell
\geq 2;
\end{eqnarray*}
  \item[{\rm (iii)}] For sup-critical distribution $\pi$ (with $\gamma=1$),
\begin{eqnarray*}
S_1 (n) &=& \frac{\|g\|_1}{n} \cdot \varphi^{-1} (n) \cdot \psi
(\log n) \cdot \Bigl[1+o (1)\Bigr], \\
S_\ell (n) &=& \frac{(\ell-2)!}{n^\ell} \cdot \varphi^{-1} (n) \cdot
\Bigl(1+o (1) \Bigr), \;  \ell \geq 2,\\
S_\ell^{(k)} (n) &=& \frac{(-1)^{k} \cdot (k+\ell-2)!}{n^{k+\ell}}
\cdot \varphi^{-1} (n) \cdot \Bigl[1+o (1)\Bigr], \; k, \ell \geq 1
\end{eqnarray*}
and
\begin{eqnarray*}
E_1 (n) &=& \|g\|_1 \cdot \varphi^{-1} (n) \cdot \psi (\log n) \cdot \Bigl[1+o (1)\Bigr], \\
E_\ell (n) &=& \frac{\varphi^{-1} (n)}{\ell \cdot (\ell-1)}\cdot
\Bigl[1+o (1)\Bigr], \; \ell \geq 2, \\
E_{\ell+} (n) &=& \frac{\varphi^{-1} (n)}{\ell-1}\cdot \Bigl( 1+o
(1) \Bigr), \; \ell \geq 2.
\end{eqnarray*}
\end{itemize}
%===========================================================================================
\end{lem}
%===========================================================================================

Let $d \geq 1$ be a fixed integer. Notice that, as $n \to \infty$,
$$
S_{\ell} (n-d) = S_\ell (n) \cdot [1+o (1)], \quad S_{\ell} (n)-S_{\ell} (n-d) = S_{\ell}^{(1)} (n_*) \cdot d
$$
for some $n_* \in (n-d, n)$, where $S_\ell^{(1)} (z)$ denotes the
derivative of $S_\ell (z)$ at $z$. Hence $S_{\ell}^{(1)} (n_*) =[1+o
(1)] \cdot S_{\ell}^{(1)} (n)$. The above lemma tells us moreover
that,
$$
S_{\ell}^{(1)} (n)=S_{\ell} (n) \cdot O (\frac{1}{n})
$$
in non-critical case or sub-critical case with $\ell \geq 1$ or
sup-critical case with $\ell \geq 2$. For sup-critical case with
$\ell=1$, we have
$$
S_1^{(1)} (n)=S_1 (n) \cdot O (\frac{1}{n \cdot \psi (\log n)})=S_1
(n) \cdot O (\frac{1}{n}).
$$
This implies the following
%===========================================================================================
\begin{lem}\label{lem: important4degree-estimate}
%===========================================================================================
For \textbf{regular} distribution $\pi$, as $n \to +\infty$ we have
$$
S_{\ell} (n-d) = S_{\ell} (n) \cdot [1+O (\frac{1}{n})],\quad
S_{\ell+1} (n) = S_{\ell} (n) \cdot O
(\frac{1}{n})
$$
for all fixed $\ell, d \geq 1$.
%===========================================================================================
\end{lem}
%===========================================================================================

%===========================================================================================
\subsection{Variation Estimation for Visiting Intensity Statistics}
%===========================================================================================
First we will estimate $\mathrm{Var\,} (R_n)$. Since $\displaystyle
R_n=\sum_x 1_{\{N_n (x) \geq 1\}}$ and $\displaystyle \Enum R_n=\sum_x [1-(1-\pi_x)^n]$,
we have
\begin{eqnarray*}
\Enum [(R_n)^2-R_n] &=& \Enum[\sum_{x \neq y} 1_{\{N_n (x) \geq 1,
N_n (y) \geq 1\}}] =\sum_{x \neq y} \Pnum (N_n (x) \geq 1, N_n (y) \geq
1)\\
&=& \sum_{x \neq y} [1-(1-\pi_x)^n-(1-\pi_y)^n-(1-\pi_x-\pi_y)^n]\\
&\leq& \sum_{x \neq y} [1-(1-\pi_x)^n][1-(1-\pi_y)^n] \leq [\Enum R_n]^2.
\end{eqnarray*}
This implies the following lemma.
%===========================================================================================
\begin{lem}\label{lem: estimate4VarR_n}
%===========================================================================================
For any distribution $\pi$, we always have $\mathrm{Var\,} (R_n) \leq \Enum R_n$.
%===========================================================================================
\end{lem}
%===========================================================================================

Now we estimate $\mathrm{Var\,} (R_{n, \, k})$. Analogously we have
\begin{eqnarray*}
&& \Enum [R_{n, \, \ell}^2 -R_{n, \, \ell}] = \sum_{x \neq y} \Pnum
(N_n (x)=\ell, N_n (y)=\ell)\\
&=& \sum_{x \neq y} \frac{n!}{(\ell!)^2 (n-2\ell)!} \cdot \pi_x^\ell
\pi_y^\ell (1-\pi_x-\pi_y)^{n-2\ell}\\
&\leq& \sum_{x \neq y} \frac{n!}{(\ell!)^2 (n-2\ell)!} \cdot
\pi_x^\ell \pi_y^\ell (1-\pi_x-\pi_y+\pi_x \pi_y)^{n-2\ell}\\
&=& \frac{n!}{(k!)^2 (n-2\ell)!} \cdot [S_\ell (n-\ell)^2-S_{2\ell} (2n-2\ell)] \leq
[1+O (\frac{1}{n})] \cdot [\frac{n^\ell \cdot S_\ell (n)}{\ell!}]^2.
\end{eqnarray*}
Noting that $\Enum R_{n, \, \ell}=C_n^\ell \cdot S_\ell
(n)=\frac{n^\ell \cdot S_\ell (n)}{\ell!} \cdot [1+O
(\frac{1}{n})]$, we have
\begin{eqnarray*}
&& \mathrm{Var\,} (R_{n, \, \ell}) = \Enum [R_{n, \, \ell}^2] -[\Enum
R_{n, \, \ell}]^2 \\
&\leq& \Enum R_{n, \, \ell} +[1+O (\frac{1}{n})] \cdot [\frac{n^\ell
\cdot S_\ell (n)}{\ell!}]^2-[1+O (\frac{1}{n})] \cdot
[\frac{n^\ell \cdot S_\ell (n)}{\ell!}]^2\\
&=& \Enum R_{n, \, \ell} +O (\frac{n^\ell \cdot S_\ell (n)}{n})
\cdot
\frac{n^\ell \cdot S_\ell (n)}{\ell!}= [1+O (\frac{\Enum R_{n, \, \ell}}{n})] \cdot \Enum R_{n, \,
\ell}.
\end{eqnarray*}
Notice that $R_{n, \, \ell} \leq R_n$ and $O (\frac{\Enum R_n}{n})=o (1)$ as $n \to \infty$, 
we always have
$$
\mathrm{Var\,} (R_{n, \, \ell})  \leq [1+o (1)] \cdot \Enum R_{n, \, \ell},
$$
But we still cannot derive an SLLN for $R_{n, \, \ell}$ directly, since
$R_{n, \, \ell}$ is not monotonic in $n$ in general.

We restate the above result as the following:
%===========================================================================================
\begin{lem}\label{lem: estimate4intensityVar}
%===========================================================================================
For \textbf{regular} distribution $\pi$, we always have
$$
\mathrm{Var\,} (R_{n, \, \ell}) \leq [1+ o(1)] \cdot \Enum R_{n, \,
\ell}
$$
as $n \to \infty$ for fixed $\ell \geq 1$.
%===========================================================================================
\end{lem}
%===========================================================================================

%===========================================================================================
\subsection{Expectation Estimation for Out-degree Statistics}
%===========================================================================================
We first estimate $\Enum \wt{R}_{n, \, k,\, \ell}$ for fixed $\ell
\geq k \geq 1$. For $\ell=1$, we have $k=1$ and $\wt{R}_{n, \, 1,\,
1}=R_{n-1, \, 1}$. We already have such estimation in the above subsections.
Thus we only need to consider the case $\ell \geq k \geq 1$ with
$\ell \geq 2$.

We introduce the following definition.
\begin{defn}
Given a vertex $x$. A finite sequence of vertices $z_1, \cdots, z_p$
of length $p \geq 2$ is called an $x$-block (of length $p$), if
$z_1=\cdots=z_{p-1}=x$ and $z_p \neq x$.
\end{defn}

For any $x$ scoring in $R_{n-1, \, \ell}$, let $\calN_n (x)$ be
the $\ell$ successive right neighbors of $x$ in the graph $G_n$; we
write them as $\tilde{x}=(x_1, \cdots, x_\ell)$; sometimes we also
regard this ordered tuple as a set: $\tilde{x} =\{x_i: i=1, \cdots, \ell\}$.
And we can partition these neighbors into the set where $x \not \in
\tilde{x}$ and the set where $x \in \tilde{x}$. We write
\begin{eqnarray*}
I_{n, \, k, \, \ell}^{(1)} (x) &:=& \sum_{\#(\tilde{x})=k, \, x \not
\in \tilde{x}} \Pnum (N_{n-1} (x)= \ell, \calN_n (x)=\tilde{x}),\\
I_{n, \, k, \, \ell}^{(2)} (x) &:=& \sum_{\#(\tilde{x})=k, \, x \in
\tilde{x}} \Pnum (N_{n-1} (x)= \ell, \calN_n (x)=\tilde{x}).
\end{eqnarray*}
Then $\Pnum (N_{n-1} (x)= \ell, D_n (x)=k) =I_{n, \, k, \, \ell}^{(1)} (x)
+I_{n, \, k, \, \ell}^{(2)} (x)$.

For the probability $I_{n, \, k, \, \ell}^{(1)} (x)$, suppose the
detailed structure of the string $\xi_1, \cdots, \xi_n$ is as
the following: the first $x$ appears at step $a_1+1$ for some $a_1
\geq 0$; after this $x$ it follows directly some vertex $x_1 \neq
x$ which gives a first contribution in the out-degree of $x$. After
the occurrence of the edge $x \to x_1$, it follows $a_2 \geq 0$
non-$x$ vertices and then an $x$ and an edge $x \to x_2$ for some
$x_2 \neq x$ and so on. Thus we get a sequence of non-negative
integers $a_1, \cdots, a_{\ell+1}$ and a sequence of vertices $x_1,
\cdots, x_\ell$. And the structure of $\xi_1, \cdots,
\xi_n$ is as the first type listed below (where $*$ denotes a
non-$x$ vertex)
\begin{equation}
*\cdots*, \hbox{edge } x \to x_1, \cdots, \hbox{edge } x
\to x_\ell, *\cdots*, \xi_n
\end{equation}
for $a_{\ell+1} \geq 1$ (and without restriction on $\xi_n$) or
as the second type listed below
\begin{equation}
*\cdots*, \hbox{edge } x \to x_1, \cdots, *\cdots*, \hbox{edge } x \to x_\ell
\end{equation}
for $a_{\ell+1}=0$ (meaning $\xi_n=x_\ell$). For the first type, we
clearly have
\begin{equation}\label{eq: solution4outdegree-1}
a_1+\cdots+a_{\ell+1}=n-2 \ell
\end{equation}
with $a_i \geq 0$ and $a_{\ell+1} \geq 1$; the number of such
integer solutions $(a_1, \cdots, a_{\ell+1})$ for equation (\ref{eq:
solution4outdegree-1}) is $C_{n-\ell-1}^\ell$. Also, the probability
of the first type for given such integer solution $(a_1, \cdots,
a_{\ell+1})$ and $\tilde{x}$ is $\pi_x^\ell \cdot (1-\pi_x)^{n-2 \ell-1} \cdot \prod_{j=1}^\ell x_j$.

For the second type, we have
\begin{equation}\label{eq: solution4outdegree-2}
a_1+\cdots+a_{\ell}=n-2\ell
\end{equation}
with $a_i \geq 0$; the number of such integer solutions $(a_1,
\cdots, a_{\ell})$ for equation (\ref{eq: solution4outdegree-2}) is
$C_{n-\ell-1}^{\ell-1}$. Also, the probability of the second type
for given such integer solution $(a_1, \cdots, a_{\ell})$ and
$\tilde{x}$ is $\pi_x^\ell \cdot (1-\pi_x)^{n-2 \ell} \cdot \prod_{j=1}^\ell x_j$.
Let's put
\begin{equation}
S_{k, \, \ell}^{(x)} :=\sum_{x \in \tilde{x}, \, \#(\tilde{x})=k}
\prod_{j=1}^\ell \pi_{x_j}.
\end{equation}
Then $I_{n, \, k, \, \ell}^{(1)} (x)$ can be formulated as
\begin{eqnarray*}
&& \Bigl[C_{n- \ell-1}^{\ell} \cdot \pi_x^\ell \cdot (1-\pi_x)^{n-2
\ell-1} + C_{n- \ell -1}^{\ell -1} \cdot \pi_x^\ell \cdot
(1-\pi_x)^{n-2 \ell} \Bigr] \cdot \sum_{x \not  \in \tilde{x}, \,
\#(\tilde{x})=k} \prod_{j=1}^\ell \pi_{x_j}\\
&=& \Bigl[C_{n- \ell-1}^{\ell} \cdot \pi_x^\ell \cdot (1-\pi_x)^{n-2
\ell-1} + C_{n- \ell -1}^{\ell -1} \cdot \pi_x^\ell \cdot
(1-\pi_x)^{n-2 \ell} \Bigr] \cdot (S_{k, \, \ell} (\pi)-S_{k, \,
\ell}^{(x)}).
\end{eqnarray*}

Note that $\displaystyle S_{k, \, \ell}^{(x)} = \sum_{r=k-1}^{\ell-1} C_\ell^r
\cdot \pi_x^{\ell-r} \cdot \sum_{x \not \in \tilde{z}, \, \#(\tilde{z})=k-1} \prod_{j=1}^r \pi_{x_j}$ is just
$$
S_{k, \, \ell}^{(x)}= \sum_{r=k-1}^{\ell-1} C_\ell^r \cdot \pi_x^{\ell-r}
\cdot [S_{k-1, r} (\pi)-S_{k-1, r}^{(x)}] \leq \sum_{r=k-1}^{\ell-1} C_\ell^r \cdot \pi_x^{\ell-r} \cdot
S_{k-1, r} (\pi).
$$
Thus if we put $\displaystyle \Delta_{n, \, k, \, \ell}^{(1)} :=\sum_x I_{n, \, k, \, \ell}^{(1)}
(x) -S_{k, \, \ell} (\pi) \cdot C_{n-\ell-1}^\ell \cdot S_{\ell} (n-\ell-1)$,
then by Lemma \ref{lem: important4degree-estimate}, $|\Delta_{n, \, k, \, \ell}^{(1)}|$ is bounded by
\begin{eqnarray*}
 & & S_{k, \, \ell} (\pi) \cdot C_{n-\ell-1}^{\ell-1} \cdot
S_{\ell} (n-\ell) +\sum_{r=1}^{\ell-k+1} S_{k-1, \, \ell-r} (\pi) \cdot C_{\ell}^r
\cdot C_{n-\ell-1}^\ell \cdot S_{\ell+r} (n-\ell+r-1)\\
&\leq& S_{k, \, \ell} (\pi) \cdot O (n^{\ell-1}) \cdot S_{\ell} (n)
\cdot [1+O (\frac{1}{n})] +\sum_{r=1}^{\ell-k+1} S_{k-1, \, \ell-r} (\pi) \cdot C_{\ell}^r
\cdot O (n^\ell) \cdot O( \frac{S_{\ell} (n)}{n^r}),
\end{eqnarray*}
which implies $|\Delta_{n, \, k, \, \ell}^{(1)}|=O (n^{\ell-1} \cdot S_{\ell} (n))$.

Similarly, in order to calculate $I_{n, \, k, \, \ell}^{(2)} (x)$,
consider the structure of the vertex sequence $\xi_1, \cdots,
\xi_n$, which may be as the following
$$
*\cdots*, x \cdots x, *\cdots*, \cdots, x \cdots x, *\cdots*.
$$
Assume $\calN_n (x)=\tilde{x}=(x_1, \cdots, x_\ell)$ where there are
elements $x_i=x$; deleting those elements $x_i=x$ from $\tilde{x}$
we obtain $\tilde{z}=(z_1, \cdots, z_p)$ for some $p \geq k-1$ and
hence there are $p$ $x$-blocks in $\xi_1, \cdots, \xi_n$ which
begin with $x$ and end with a unique non-$x$ vertex. And the
detailed structure of $\xi_1, \cdots, \xi_n$ is as the following:
after $a_1 \geq 0$ many non-$x$ vertices, the first $x$-block
appears which contains first $b_1 \geq 1$ many $x$ and then a unique
non-$x$ vertex. After this first $x$-block, it follows $a_2 \geq 0$
many non-$x$ vertices and then the second $x$-blocks and so on.
Clearly $a_1+\cdots+a_p+a_{p+1}=n-p-\ell$ and
$b_1+\cdots+b_p=\ell$ or $ b_1+\cdots+b_p=\ell+1$; the last case
$b_1+\cdots+b_p=\ell+1$ corresponds to $\xi_n=x$, $a_{p+1}=0, b_p
\geq 2$ and $\tilde{z}=(z_1, \cdots, z_{p-1})$ (requiring
$p-1 \geq k-1=\#(\tilde{z})$). So
\begin{eqnarray*}
&& I_{n, \, k, \, \ell}^{(2)} (x) =
\sum_{\stackrel{\#(\tilde{x})=k}{x \in \tilde{x}}} \Pnum (N_{n-1}
(x)=\ell, \calN_n (x)=\tilde{x})\\
&=& \sum_{p=k-1}^{\ell-1} \sum_{\stackrel{\#(\tilde{z})=k-1}{z_i
\neq x}} \; \sum_{\stackrel{b_1+\cdots+b_p=\ell}{b_i \geq 1}} \;
\sum_{\stackrel{a_1+\cdots+a_{p+1}=n-p-\ell}{a_i \geq
0}} \pi_x^\ell (\prod_{j=1}^p \pi_{z_j}) (1-\pi_x)^{a_1+\cdots+a_p+(a_{p+1}-1)^+} \\
&+& \sum_{p=k}^{\ell-1} \sum_{\stackrel{\#(\tilde{z})=k-1}{z_i \neq
x}} \; \sum_{\stackrel{b_1+\cdots+b_p=\ell+1}{b_i \geq 1, \, b_p
\geq 2}} \; \sum_{\stackrel{a_1+\cdots+a_{p}=n-p-\ell}{a_i \geq
0}} \pi_x^{\ell+1} (\prod_{j=1}^{p-1} \pi_{z_j}) (1-\pi_x)^{a_1+\cdots+a_p}\\
&=& \sum_{p=k-1}^{\ell-1} C_{\ell-1}^{p-1} \cdot C_{n-\ell-1}^p \cdot (1-\pi_x)^{n-p-\ell-1}
\cdot \pi_x^\ell \cdot [S_{k-1, \, p}
(\pi)-S_{k-1, \, p}^{(x)}]\\
 &+& \sum_{p=k-1}^{\ell-1} C_{\ell-1}^{p-1} \cdot C_{n-\ell-1}^{p-1} \cdot (1-\pi_x)^{n-p-\ell}
\cdot \pi_x^\ell \cdot [S_{k-1, \, p}
(\pi)-S_{k-1, \, p}^{(x)}]\\
&+& \sum_{p=k}^{\ell-1} C_{\ell-1}^{p-1} \cdot C_{n-\ell-1}^{p-1}
\cdot (1-\pi_x)^{n-p-\ell} \cdot \pi_x^{\ell+1} \cdot [S_{k-1, \, p}
(\pi)-S_{k-1, \, p}^{(x)}].
\end{eqnarray*}
Therefore $\sum_x I_{n, \, k, \, \ell}^{(2)} (x)$ is bounded by
$$
C \cdot \left\{ \sum_{p=k-1}^{\ell-1}  [n^p \cdot S_{\ell} (n-p-1)+n^{p-1} \cdot S_{\ell} (n-p)]
+\sum_{p=k}^{\ell-1} n^{p-1} \cdot S_{\ell+1} (n-p+1) \right\}
$$
with some constant $C$. Now one can prove $\displaystyle \sum_x
I_{n, \, k, \, \ell}^{(2)} (x)=O (\cdot n^{\ell-1} \cdot S_{\ell} (n))$ by noting our Lemma 
\ref{lem: important4degree-estimate}. Therefore we have the following
%===========================================================================================
\begin{lem}\label{lem: estimate4outdegree-1}
%===========================================================================================
For regular distribution $\pi$, we always have
\begin{equation}
\Enum \wt{R}_{n, \, k, \, \ell}= S_{k, \, \ell} (\pi) \cdot
\frac{n^\ell \cdot S_{\ell} (n)}{\ell!} \cdot [1+O (\frac{1}{n})]
\end{equation}
as $n \to \infty$ for $1 \leq k \leq \ell$ and $\ell \geq 2$ fixed.
%===========================================================================================
\end{lem}
%===========================================================================================

%===========================================================================================
\subsection{Variation Estimation for Out-degree Statistics}
%===========================================================================================
Now we study the variation of $\wt{R}_{n, \, k, \, \ell}$ for fixed $1
\leq k \leq \ell$ and large enough $n$. As before, we have
\begin{eqnarray*}
\Enum (\wt{R}_{n, \, k, \, \ell}^2-\wt{R}_{n, \, k, \, \ell}) &=&
\sum_{(x,\,y): x \neq y} \Pnum \Bigl( N_{n-1} (x)=\ell, D_n
(x)=k, N_{n-1} (y)=\ell, D_n (y)=k \Bigr)\\
&=:& \sum_{(x,\,y): x \neq y} J_{n, \, k, \, \ell} (x, y).
\end{eqnarray*}

We introduce the following definition.
\begin{defn}
Given two distinct vertices $x \neq y$. A finite sequence of
vertices $z_1, \cdots, z_p$ of length $p \geq 2$ is called an
$(x,y)$-block (of length $p$), if $z_1, \cdots, z_{p-1} \in \{x,
y\}$ and $z_p \not \in \{x, y\}$.
\end{defn}

We now calculate the probability $J_{n, \, k, \, \ell} (x, y)$ for
fixed $1 \leq k \leq \ell$ and distinct vertices $x \neq y$. In the
calculations below, we will denote by
$$
\tilde{x}=(x_1, \cdots, x_{\ell})=\calN_n (x), \quad \tilde{y}=(y_1,
\cdots, y_{\ell})=\calN_n (y)
$$
the right neighbors of $x$ and $y$ respectively in the graph $G_n$
with the restriction
\begin{equation}\label{eq: out-degree-restriction}
\#\{x_1, \cdots, x_{\ell}\}=\#\{y_1, \cdots, y_{\ell}\}=k.
\end{equation}
Clearly, such neighbors $\tilde{x}$ and $\tilde{y}$ can be
partitioned into two sets: the first set is such that both
$\tilde{x}$ and $\tilde{y}$ have no elements being $x$ or $y$; the
second set is such that there is some element either of $\tilde{x}$
or of $\tilde{y}$ being $x$ or $y$. Corresponding to such a
partition, we write
\begin{eqnarray*}
J^{(1)}_{n,\, k,\, \ell} (x,\, y) &:=& \sum_{\hbox{each } x_i,\, y_i
\not \in \{x, \, y\}} \Pnum \Bigl( N_{n-1} (x)=\ell, \calN_n
(x)=\tilde{x},
N_{n-1} (y)=\ell, \calN_n (y)=\tilde{y} \Bigr)\\
J^{(2)}_{n,\, k,\, \ell} (x,\, y) &:=& \sum_{\exists i, \, x_i
\hbox{ or } y_i \in \{x, \, y\}} \Pnum \Bigl( N_{n-1} (x)=\ell,
\calN_n (x)=\tilde{x}, N_{n-1} (y)=\ell, \calN_n (y)=\tilde{y}
\Bigr),
\end{eqnarray*}
where the restriction (\ref{eq: out-degree-restriction}) is omitted
in the summations for simplicity of presentation. Then $J_{n, \, k, \, \ell} (x, y)) = J^{(1)}_{n,\, k,\, \ell} (x,\, y)+
J^{(2)}_{n,\, k,\, \ell} (x,\, y)$.

We calculate $J^{(1)}_{n,\, k,\, \ell} (x,\, y)$ first. Put
$\displaystyle S_{k,\, \ell}^{(x,y)} :=\sum_{\tilde{x}: \exists x_i \in \{x, \,
y\}} \prod_{j=1}^\ell x_j$. Considering a typical realization of the random sequence $\xi_1,
\cdots, \xi_n$ satisfying the obvious restriction in the calculation
of $J^{(1)}_{n,\, k,\, \ell} (x,\, y)$, where both $x$ and $y$ show
up exactly $\ell$-times in the first $n-1$ steps, leaving their
neighbors $\tilde{x}, \tilde{y}$. Since the neighbors $\tilde{x},
\tilde{y}$ have no elements being $x$ or $y$, the random sequence
$\xi_1, \cdots, \xi_n$ have exactly $\ell$ $x$-blocks and $\ell$
$y$-blocks. Noting that the relative disposition of these $\ell$
$x$-blocks (respectively, $y$-blocks) is uniquely determined by
$\tilde{x}$ (respectively, $\tilde{y}$), there are $C_{2 \ell}^\ell$
kinds of relative dispositions of these $\ell$ $x$-blocks and $\ell$
$y$-blocks. When the relative disposition of of these $\ell$
$x$-blocks and $\ell$ $y$-blocks is fixed, we say that we have $2
\ell$ $(x,y)$-blocks in the random sequence. And the detailed
structure of the random sequence is of the following two type: (1)
after $a_1 \geq 0$ many non-$(x,y)$ vertices, the first
$(x,y)$-block appears, and then follows $a_2 \geq 0$ many
non-$(x,y)$ vertices and so on; and after the last $(x,y)$-block, it
follows $a_{2\ell+1} \geq 1$ many non-$(x,y)$ vertices; i.e., the
structure is as the following (where $*$ denotes non-$(x,y)$ vertex)
\begin{equation}
*\cdots*, (x,y)-\hbox{block}, *\cdots*, \cdots, (x,y)-\hbox{block}, *\cdots*
\end{equation}
(2) the detailed structure is almost the same as the first type (1)
with the only modification that $a_{2\ell+1}=0$, i.e., the structure
is as the following (where $*$ denotes non-$(x,y)$ vertex)
\begin{equation}
*\cdots*, (x,y)-\hbox{block}, *\cdots*, \cdots, (x,y)-\hbox{block}
\end{equation}
It is clear that there are $C_{n-2\ell-1}^{2 \ell}$ solutions to the
equation
$$
a_1+ \cdots+a_{2 \ell+1}=n-4\ell
$$
with $a_i \geq 0$ and $a_{2\ell+1} \geq 1$. Also, the probability of
the random sequence being the first type structure is $\pi_x^\ell \cdot \pi_y^\ell \cdot (1-\pi_x-\pi_y)^{n-4\ell-1} \cdot
\prod_{j=1}^\ell \pi_{x_j} \cdot \prod_{j=1}^\ell \pi_{y_j}$.

Similarly, there are $C_{n-2\ell-1}^{2 \ell-1}$ solutions to the
equation
$$
a_1+ \cdots+a_{2 \ell}=n-4\ell
$$
with $a_i \geq 0$. And the probability of the random sequence being
the second type structure is $\pi_x^\ell \cdot \pi_y^\ell \cdot (1-\pi_x-\pi_y)^{n-4\ell} \cdot
\prod_{j=1}^\ell \pi_{x_j} \cdot \prod_{j=1}^\ell \pi_{y_j}$.

Summing up, we have
\begin{eqnarray*}
&& J^{(1)}_{n,\, k,\, \ell} (x,\, y) = C_{2 \ell}^\ell \cdot
\Bigl[C_{n-2\ell-1}^{2 \ell} \cdot \pi_x^\ell \cdot \pi_y^\ell \cdot
(1-\pi_x-\pi_y)^{n-4\ell-1} \\
&&+C_{n-2\ell-1}^{2 \ell-1} \cdot \pi_x^\ell \cdot \pi_y^\ell \cdot
(1-\pi_x-\pi_y)^{n-4\ell}\Bigr] \cdot
\sum_{x_i, y_j \not\in \{x,y\}} \prod_{j=1}^\ell \pi_{x_j} \cdot \prod_{j=1}^\ell \pi_{y_j}
\end{eqnarray*}
which is
\begin{eqnarray*}
J^{(1)}_{n,\, k,\, \ell} (x,\, y) &=& C_{2\ell}^\ell \cdot C_{n-2\ell}^{2\ell} \cdot
(1-\pi_x-\pi_y)^{n-4\ell-1} \cdot (\pi_x \pi_y)^\ell \cdot [S_{k,
\ell} (\pi) -S_{k,\, \ell}^{(x,y)}]^2\\
&+& C_{2\ell}^\ell \cdot C_{n-2\ell-1}^{2\ell-1} \cdot
(1-\pi_x-\pi_y)^{n-4\ell} \cdot (\pi_x \pi_y)^\ell \cdot [S_{k,
\ell} (\pi) -S_{k,\, \ell}^{(x,y)}]^2\\
&\leq& C_{2\ell}^\ell \cdot C_{n-2\ell}^{2\ell} \cdot
[(1-\pi_x)^{n-4\ell-1} \cdot (\pi_x)^\ell] \cdot
[(1-\pi_y)^{n-4\ell-1} \cdot (\pi_y)^\ell] \cdot [S_{k, \ell}
(\pi)]^2\\
&+& C_{2\ell}^\ell \cdot C_{n-2\ell-1}^{2\ell-1} \cdot
[(1-\pi_x)^{n-4\ell} \cdot (\pi_x)^\ell] \cdot [(1-\pi_y)^{n-4\ell}
\cdot (\pi_y)^\ell] \cdot [S_{k, \ell} (\pi)]^2.
\end{eqnarray*}
Hence $\displaystyle J^{(1)}_{n,\, k,\, \ell} :=\sum_{(x,y): x \neq y} J^{(1)}_{n,\,
k,\, \ell} (x,\, y)$ is bounded by
$$
J^{(1)}_{n,\, k,\, \ell} \leq [S_{k, \ell} (\pi)]^2 \cdot C_{2\ell}^\ell \cdot \Bigl[
C_{n-2\ell-1}^{2\ell-1} \cdot S_{\ell} (n-3\ell)^2+
C_{n-2\ell}^{2\ell} \cdot S_{\ell} (n-3\ell-1)^2\Bigr].
$$
For $\ell \geq 2$, we have $\displaystyle J^{(1)}_{n,\, k,\, \ell} \leq
[S_{k, \ell} (\pi)]^2 \cdot \Bigl[\frac{n^\ell \cdot S_{\ell} (n)}{\ell!}
\Bigr]^2 \cdot \Bigl[1+O(\frac{1}{n})\Bigr]$.

Similarly, in calculating $\displaystyle \Pnum \Bigl( N_{n-1} (x)=\ell, \calN_n (x)=\tilde{x}, N_{n-1}
(y)=\ell, \calN_n (y)=\tilde{y}\Bigr)$ with the condition $\exists i, \, x_i \hbox{ or } y_i \in \{x, \,
y\}$, first consider the situations $k=1$ and $k=2$.

If $k=1$, then the detailed structure of $\xi_1, \cdots, \xi_n$ can
be one of the following forms (where $*$ denotes a non-$(x,y)$
vertex):
\begin{eqnarray}
\label{eq: structure4(k=1)-1} *\cdots*, y-\hbox{block}, *\cdots*,
y-\hbox{block}, *\cdots*, x \cdots x\\
\label{eq: structure4(k=1)-2} *\cdots*, x-\hbox{block}, *\cdots*,
x-\hbox{block}, *\cdots*, y \cdots y\\
\label{eq: structure4(k=1)-3} *\cdots*, (x, y)-\hbox{block},
*\cdots*, (x, y)-\hbox{block}, *\cdots*, x, y\\
\label{eq: structure4(k=1)-4} *\cdots*, (x, y)-\hbox{block},
*\cdots*, (x, y)-\hbox{block}, *\cdots*, y, x.
\end{eqnarray}
And according to the above formulations, we have
\begin{eqnarray*}
J^{(2)}_{n,\, 1,\, \ell} (x, y) &=& C_{n-2\ell-1}^\ell \cdot
\pi_x^{\ell+1} \cdot \pi_y^\ell \cdot (1-\pi_x-\pi_y)^{n-3 \ell-1}
\cdot [S_{1, \ell} (\pi) -S_{1, \ell}^{(x,y)}]\\
&+& C_{n-2\ell-1}^\ell \cdot \pi_x^{\ell} \cdot \pi_y^{\ell+1} \cdot
(1-\pi_x-\pi_y)^{n-3 \ell-1} \cdot [S_{1, \ell} (\pi) -S_{1, \ell}^{(x,y)}]\\
&+& C_\ell^1 \cdot C_{n-2\ell-1}^\ell \cdot \pi_x^{\ell} \cdot
\pi_y^{\ell+1} \cdot (1-\pi_x-\pi_y)^{n-3 \ell-1} \cdot [S_{1, \ell}
(\pi) -S_{1, \ell}^{(x,y)}]\\
&+& C_\ell^1 \cdot C_{n-2\ell-1}^\ell \cdot \pi_x^{\ell+1} \cdot
\pi_y^{\ell} \cdot (1-\pi_x-\pi_y)^{n-3 \ell-1} \cdot [S_{1, \ell}
(\pi) -S_{1, \ell}^{(x,y)}]
\end{eqnarray*}
which implies $\displaystyle J^{(2)}_{n,\, 1,\, \ell} :=\sum_{x \neq y} J^{(2)}_{n,\, 1,\,
\ell} (x, y)$ is bounded by
\begin{eqnarray*}
J^{(2)}_{n,\, 1,\, \ell} &\leq& 2 \cdot S_{1, \ell} (\pi) \cdot (\ell+1) \cdot
C_{n-2\ell-1}^\ell \cdot S_{\ell+1} (n-2\ell) \cdot S_{\ell}
(n-2\ell-1)\\
&=& O (n^{\ell-1} \cdot S_{\ell} (n)^2)=O (\frac{1}{n^{\ell+1}})
\cdot [\Enum \wt{R}_{n, \, 1, \, \ell}]^2=o (\frac{1}{n^\ell}) \cdot
\Enum \wt{R}_{n, \, 1, \, \ell}
\end{eqnarray*}
since we have proved $\displaystyle \Enum \wt{R}_{n, \, 1, \, \ell}
=S_{1, \, \ell} (\pi) \cdot \frac{n^\ell \cdot S_\ell (n)}{\ell!} \cdot
[1+O (\frac{1}{n})]$ by our Lemma \ref{lem: estimate4outdegree-1}; Therefore
$\displaystyle \mathrm{Var\;} (\wt{R}_{n, \, 1, \, \ell}) \leq [1+o (1)] \cdot
\Enum \wt{R}_{n, \, 1, \, \ell}$ as $n \to \infty$ for fixed $\ell \geq 2$.

For $k=2$, $J^{(2)}_{n,\, 2,\, \ell} (x, y)$ ($x \neq y$) can be calculated
in the following way. There are 4 patterns for the edges
with starting vertex $x$ (restricted to $D_n (x)=2$),
where $*$ denotes some non-$(x,y)$-vertex: (1) $x \to *, x \to *$;
(2) $x \to x, x \to *$; (3) $x \to y, x \to *$; (4) $x \to x, x \to
y$; A similar result holds for the vertex $y$. Thus there are $4
\times 4-1=15$ patterns for the edges with starting vertices $x,y$
(restricted to the conditions $D_n (x)=2, D_n (y)=2$ and that there
exists some right neighbor of $x$ or $y$ being exactly $x$ or $y$).
For each of these 15 patterns, there are in general $p \leq 2 \cdot
\ell-1$ $(x,y)$-blocks in the random sequence $\xi_1, \cdots,
\xi_n$, which results into a factor $C_{n-2\ell-1}^p=O (n^p) \leq O
(n^{2 \cdot \ell-1})$ along with the related probability bounded by
$\pi_x^\ell \cdot \pi_y^\ell \cdot (1-\pi_x-\pi_y)^{n-4 \ell}$,
therefore $\displaystyle J^{(2)}_{n,\, 2,\, \ell} (x, y) \leq
O (n^{2 \cdot \ell-1}) \cdot \pi_x^\ell \cdot \pi_y^\ell \cdot
(1-\pi_x-\pi_y)^{n-4 \ell}$ and hence
\begin{eqnarray*}
J^{(2)}_{n,\, 2,\, \ell} &:=& \sum_{x \neq y} J^{(2)}_{n,\, 2,\,
\ell} (x, y) \leq \sum_{x \neq y} O (n^{2 \cdot \ell-1}) \cdot
\pi_x^\ell \cdot \pi_y^\ell \cdot (1-\pi_x-\pi_y)^{n-4 \ell}\\
&\leq& \sum_{x \neq y} O (n^{2 \cdot \ell-1}) \cdot
\pi_x^\ell \cdot \pi_y^\ell \cdot (1-\pi_x)^{n-4 \ell} \cdot (1-\pi_y)^{n-4 \ell}\\
 &\leq& O (n^{2 \cdot \ell-1} \cdot S_\ell (n-3 \ell)^2)=o (1) \cdot
 \Enum \wt{R}_{n,\, 2,\, \ell}.
\end{eqnarray*}
Thus $\mathrm{Var\;} (\wt{R}_{n, \, 2, \, \ell}) \leq [1+o (1)] \cdot
\Enum \wt{R}_{n, \, 2, \, \ell}$ as $n \to \infty$ for fixed $\ell \geq 2$.

For $k \geq 3$, we first partition the case into two cases: (A) $x_\ell, y_\ell \not\in \{x,y\}$; (B) either
$x_\ell \in \{x,y\}$ or $y_\ell \in \{x,y\}$ (this means $\xi_n \in
\{x, y\}$). For the case (A), in the realization of the random
sequence $\xi_1, \cdots, \xi_n$, suppose there are $p$ many
$(x,y)$-blocks; clearly $p \leq 2\ell-1$. For the case (B), in the
realization of the random sequence $\xi_1, \cdots, \xi_n$, suppose
there are $p$ many $(x,y)$-blocks; clearly $p \leq 2 \ell-1$.
Similar to the calculation of $J^{(2)}_{n,\, 2,\, \ell}$, we would
have
$$
J^{(2)}_{n,\, k,\, \ell} (x, y) \leq O (n^{2 \ell-1}) \cdot
\pi_x^{\ell} \cdot \pi_y^{\ell} \cdot (1-\pi_x-\pi_y)^{n-4 \ell}.
$$
Hence $\displaystyle J^{(2)}_{n,\, k,\, \ell} := \sum_{(x,y): x \neq y} J^{(2)}_{n,\,
k,\, \ell} (x, y) \leq [n^\ell \cdot S_{\ell} (n)]^2 \cdot O (\frac{1}{n})=o (1)
\cdot \Enum \wt{R}_{n ,\, k, \,\ell}$.

Summing up the above results, we have
\begin{lem}\label{lem: estimate4outdegreeVar}
For \textbf{regular} distribution $\pi$, we always have
$$
\mathrm{Var\;} (\wt{R}_{n, \, k, \, \ell}) \leq \Bigl[1+o (1)\Bigr] \cdot \Enum \wt{R}_{n ,\,
k, \,\ell}
$$
as $n \to \infty$ for fixed $1 \leq k \leq \ell$.
\end{lem}
%===========================================================================================
\section{Proofs for the Main Theorems \ref{thm: SLLN4Rn}--\ref{thm: main-sub-critical}}\label{sec:4}
%===========================================================================================
Now we shall derive the detailed range-renewal structure for i.i.d. sequences. For simplicity
of presentation, the proof for the results of the induced undirected graph $\wh{G}_n$ is omitted
here since it is rather similar to the directed graph case. The main idea of these proofs is to
exploit Lemma \ref{lem: SLLN-approach} to build a sequence of SLLNs for regular distributions.

First, for any distribution $\pi$, we always have $\mathrm{Var\,}
(R_n) \leq \Enum R_n$, which implies $\displaystyle \lim_{n \to +\infty} \frac{R_n}{\Enum R_n}=1$
almost surely in view of Lemma \ref{lem: SLLN-approach}.

Now in view of Cauchy's inequality and the estimations in Section
\ref{sec:3}, we have for any $\ell \geq 2$
\begin{eqnarray*}
\mathrm{Var\,} (R_{n, \, \ell+}) &=& \mathrm{Var\,} (R_n
-\sum_{k=1}^{\ell-1} R_{n, \, k}) \leq \ell \cdot \Bigl[ \mathrm{Var\,} (R_n) +\sum_{k=1}^{\ell-1}
\mathrm{Var\,} (R_{n, \, k}) \Bigr]\\
&\leq& \ell \cdot \Bigl[ \Enum R_n +(1+o(1)) \cdot
\sum_{k=1}^{\ell-1}
\Enum R_{n, \, k}\Bigr] \leq 2 \ell (1+o(1)) \cdot \Enum R_n,
\end{eqnarray*}
i.e.,
\begin{equation}\label{eq: base4SLLN}
\mathrm{Var\,} (R_{n, \, k+}) \leq C_1 \cdot \Enum R_n
\end{equation}
for some positive constant $C_1$. The above estimation is the
starting point of this part. Note also that, now $R_{n, \, k+}$ is
increasing in $n$ for any fixed $k \geq 2$; so it's possible to
obtain SLLNs for such sequences by our Lemma
\ref{lem: SLLN-approach}.

%===========================================================================================
\subsection{Non-Critical Case: $\gamma=\gamma (\pi) \in (0, 1)$}\label{sec:4.1}
%===========================================================================================
For the non-critical case, we have already known that for $\ell \geq
2$
$$
\Enum R_{n, \, \ell+} =\frac{\Gamma (\ell-\gamma)}{(\ell-1)!} \cdot
\varphi^{-1} (n) \cdot \Bigl( 1+o(1) \Bigr)
$$
and $\Enum R_n=\Gamma (1-\gamma) \cdot \varphi^{-1} (n) \cdot \Bigl(
1+o(1) \Bigr)$. Hence by (\ref{eq: base4SLLN}) we obtain
$$
\mathrm{Var\,} (R_{n, \, \ell+}) \leq C \cdot \Enum R_{n, \, \ell+}
$$
for some constant $C>0$ (which may depend on $\ell$). Now noting
$$
\frac{\Enum R_{n, \, \ell}}{\Enum R_n}=\frac{\gamma \cdot \prod
\limits_{j=1}^{\ell-1} (j-\gamma)}{\ell!} \cdot \Bigl( 1+o(1)
\Bigr),
$$
we can easily derive $\displaystyle \lim_{n \to \infty} \frac{R_{n, \ell}}{R_n}=r_{\ell} (\gamma), \;
\ell \geq 1$ almost surely.

In order to obtain information for $\wt{R}_{n, \, k, \, \ell}$ with
$1 \leq k \leq \ell$, we put
$$
\wt{R}_{n, \, k+, \, \ell+} :=\sum_{i \geq k \hbox{ or } j \geq
\ell} \wt{R}_{n, \, i, \, j}.
$$
It's clear that $\wt{R}_{n, \, k+, \, \ell+}$ is
non-decreasing in $n$ for fixed $k, \ell$. Also, noting
$$
\wt{R}_{n, \, k+, \, \ell+}=R_{n-1}-\sum_{i<k \hbox{ and } j <\ell}
\wt{R}_{n, \, i, \, j},
$$
where the right hand side has only $N=N (k, \, \ell)<\infty$ terms,
we have
\begin{eqnarray*}
\mathrm{Var\;} (\wt{R}_{n, \, k+, \, \ell+}) &\leq& N^2 \cdot \Bigl[
\mathrm{Var\;}  (R_{n-1})+\sum_{i<k \hbox{ and } j <\ell}
\mathrm{Var\;} (\wt{R}_{n, \, i, \, j}) \Bigr]\\
&\leq& N^2 \cdot \Bigl[ \Enum (R_{n-1})+\sum_{i<k \hbox{ and } j
<\ell} \Bigl[1+o (1)\Bigr] \Enum \wt{R}_{n, \, i, \, j} \Bigr]\\
&\leq& 2 \cdot N^2 \cdot \Bigl[1+o (1)\Bigr] \cdot \Enum R_{n-1}
\end{eqnarray*}
as $n \to \infty$; it is easy to see that $\Enum R_{n-1}$ and $\Enum
\wt{R}_{n, \, k+, \, \ell+}$ are of the same order. Hence we have
SLLN for $\wt{R}_{n, \, k+, \, \ell+}$: almost surely $\displaystyle
\lim_{n \to \infty} \frac{\wt{R}_{n, \, k+, \, \ell+}}{\Enum
\wt{R}_{n, \, k+, \, \ell+}}=1$. From this we can easily obtain Theorem
\ref{thm: main-non-critical} since we have
\begin{eqnarray}
R_{n, \, \ell} &=& R_{n, \, \ell+} -R_{n, \, (\ell+1)+}, \\
\nonumber \wt{R}_{n, \, k, \, \ell} &=&\Bigl[ \wt{R}_{n, \, k+, \,
(\ell+1)+} -\wt{R}_{n, \, (k+1)+, \, (\ell+1)+} \Bigr]\\
&-& \Bigl[ \wt{R}_{n, \, k+, \, \ell+} -\wt{R}_{n, \, (k+1)+, \,
\ell+} \Bigr],
\end{eqnarray}
where by convention $\wt{R}_{n, \, k+, \,
\ell+}=R_{n-1, \ell+}$ if $k>\ell$ and $R_{n, 1+}=R_n$.
%===========================================================================================
\subsection{Sub-Critical Case: $\gamma (\pi)=0$}\label{sec:4.2}
%===========================================================================================
In the same spirit as above, we can prove the following SLLNs
$$
\lim_{n \to \infty} \frac{R_{n, \, \ell+}}{\Enum R_{n, \, \ell+}}
= 1 \hbox{ and }
\lim_{n \to \infty} \frac{\wt{R}_{n, \, k+, \, \ell+}}{\Enum
\wt{R}_{n, \, k+, \, \ell+}} = 1.
$$
A careful calculation reveals that
$$
\Enum R_{n, \, \ell+}=[1+o (1)] \cdot \Enum R_n \hbox{ and } \Enum \wt{R}_{n,
\, k+, \, \ell+} =[1+o (1)] \cdot \Enum R_n,
$$
which implies $\displaystyle \lim_{n \to \infty} \frac{R_{n, \, \ell+}}{R_n} = 1$and $\displaystyle
\lim_{n \to \infty} \frac{\wt{R}_{n, \, k+, \, \ell+}}{R_n} = 1$. Hence the results in Theorem
\ref{thm: main-sub-critical} hold true.
%===========================================================================================
\subsection{Sup-Critical Case: $\gamma (\pi)=1$}\label{sec:4.3}
%===========================================================================================
For this case, we have (\ref{eq: sup-critical-order}) (see Corollary
\ref{cor: sup-critical-order}). So it's clear that
$$
\lim_{n \to \infty} \frac{\log \Enum R_n}{\log n} = 1, \quad \lim_{n \to \infty}
\frac{\log \Enum R_{n, \, \ell+}}{\log n} = 1, \quad \ell \geq 2
$$
since $\Enum R_n=O (\varphi^{-1} (n) \cdot \psi (\log n)), \Enum
R_{n, \, \ell+}=O (\varphi^{-1} (n))$ for $\ell \geq 2$ as $n \to
\infty$. Therefore we derive easily that $\mathrm{Var\,} (R_{n, \, \ell+}) \leq C \cdot
\Bigl( \Enum R_{n, \, \ell+} \Bigr)^{3/2}$ for some constant $C>0$ (which may depend on
$\ell \geq 2$), which yields an SLLN for $R_{n, \, \ell+}$ by our Lemma
\ref{lem: SLLN-approach}. In the same spirit, we have also an SLLN
for $\wt{R}_{n, \, k+, \, \ell+}$.

Now a detailed calculation of $\Enum R_{n, \, \ell+}$ and $\Enum
\wt{R}_{n, \, k+, \, \ell+}$ reveals the results in Theorem
\ref{thm: main-sup-critical}.

%\section*{Acknowledgements} The
%first author would like to thank Prof. De-Jun Feng for helpful
%discussion and comments. He also expresses the gratitude to Prof.
%Derriennic for sending him a copy of ref. \cite{Derriennic}.
%========================================================================================
\noindent{\sl \textbf{Acknowledgements} \quad} The second author
would like to thank Prof. De-Jun Feng, Prof. Narn-Rueih Shieh for helpful discussions during
his short visit at CUHK in 2012; He also thanks his colleagues Prof. Yun-Xin Zhang and Prof.
Yi-Jun Yao for helpful discussions and comments; He would like to expresses the gratitude to
Prof. Derriennic for sending him a copy of \cite{Derriennic}. This work is in part supported by
NSFC (No. 11001173, No. 11271255, No. 10701026 and No. 11271077) and the Laboratory of
Mathematics for Nonlinear Science, Fudan University.



%% References


\bibliographystyle{elsarticle-num}

\begin{thebibliography}{00}


\bibitem{AH00} Adamic, Lada A.; Huberman, Bernardo A.;
{\em Power-Law Distribution of the World Wide Web}, Science
\textbf{287}, 2115 (2000).

\bibitem{AB02} Albert, R.; Barab\'{a}si, A.-L.;
{\em Statistical mechanics of complex networks}, Rev. Mod. Phys.
\textbf{74} (2002), pp. 47--97.

\bibitem{Athreya85} Athreya, K. B.;
{\em On the Range of Recurrent Markov Chains}, Statist. Probab.
Lett. \textbf{3} (1985), no. 3, pp. 143--145. MR0801860

\bibitem{B09} Barab\'{a}si, A.-L.;
{\em Scale-Free Networks: A Decade and Beyond}, Science
\textbf{325}, 412 (2009).

\bibitem{BA99} Barab\'{a}si, A.-L.; Albert, R.;
{\em Emergence of Scaling in Random Networks}, Science \textbf{286},
509 (1999).

\bibitem{BK02} Bass, Richard F.; Kumagai, Takashi;
{\em Laws of the iterated logarithm for the range of random walks in
two and three dimensions}, Ann. Probab. \textbf{30} (2002), no. 3,
1369--1396. MR1920111

\bibitem{CI78} Chosid, Leo; Isaac, Richard;
{\em On the Range of Recurrent Markov Chains}, Ann. Probab.
\textbf{6} (1978), no. 4, pp. 680--687. MR0474507

\bibitem{CI80} Chosid, Leo; Isaac, Richard;
{\em Correction to: ``On the range of recurrent Markov chains''
[Ann. Probab. \textbf{6} (1978), no. 4, 680--687; MR 57 \#14146]}.
Ann. Probab. \textbf{8} (1980), no. 5, pp. 1000. MR0600347

\bibitem{Derriennic} Derriennic, Yves;
{\em Quelques applications du th\'{e}or\`{e}me ergodique
sous-additif}. (French. English summary) Conference on Random Walks
(Kleebach, 1979) (French), pp. 183--201, 4, Ast\'{e}risque,
\textbf{74}, Soc. Math. France, Paris, 1980. MR0588163

\bibitem{DE50} Dvoretzky, A.; Erd\"{o}s, P.
{\em Some problems on random walk in space}. Proceedings of the
Second Berkeley Symposium on Mathematical Statistics and
Probability, 1950. pp. 353--367. University of California Press,
Berkeley and Los Angeles, 1951. MR0047272

\bibitem{ET} Erd\"{o}s, P.; Taylor, S. J.;
{\em Some problems concerning the structure of random walk paths},
Acta Math. Acad. Sci. Hungar. \textbf{11} (1960), 137--162. MR0121870

\bibitem{Feller} Feller, W.:
{\it An Introduction to Probability Theory and Its Applications},
2nd edition, Vol. 2. ISBN 0-471-25709-5, Wiley Publishing, Inc. (Chinese translation
edition, Posts \& Telecom Press, 2008.)

%\bibitem{GHZ76} Geman, Donald; Horowitz, Joseph; Zinn, Joel;
%{\em Recurrence of stationary sequences}, Ann. Probability
%\textbf{4} (1976), no. 3, pp. 372--381.

\bibitem{Glynn85} Glynn, Peter W.;
{\em On the Range of a Regenerative Sequence}, Stochastic Process.
Appl. \textbf{20} (1985), no. 1, pp. 105--113. MR0805118

\bibitem{GS92} Glynn, Peter; Sigman, Karl;
{\em Uniform Ces$\grave{a}$ro Limit Theorems for Synchronous
Processes with Applications to Queues}, Stochastic Process. Appl.
\textbf{40} (1992), no. 1, pp. 29--43. MR1145457

\bibitem{Hamana97} Hamana, Y.;
{\em The fluctuation result for the multiple point range of
two-dimensional recurrent random walks}, Ann. Probab. \textbf{25}
(1997), 598--639. MR1434120

\bibitem{Hamana98} Hamana, Y.;
{\em An almost sure invariance principle for the range of random
walks}, Stochastic Process. Appl. \textbf{78} (1998), 131--143.
MR1657371

\bibitem{I-K} Iosifescu, M.; Kraaikamp, C.:
{\it  Metrical Theory of Continued Fractions}, Mathematics and its
Applications, 547. Kluwer Academic Publishers, Dordrecht, 2002. MR1960327

\bibitem{JP71} Jain, N. C.; Pruitt, W. E.;
{\em The range of transient random walk}, J. Anal. Math. \textbf{24}
(1971), 369--393. MR0283890

\bibitem{JP72} Jain, N. C.; Pruitt, W. E.;
{\em The law of the iterated logarithm for the range of random
walk}, Ann. Math. Statist. \textbf{43} (1972), 1692--1697. MR0345216

\bibitem{JP72'} Jain, N. C.; Pruitt, W. E.;
{\em The range of random walk}, Proc. Sixth Berkeley Symp. Math.
Statist. Probab. \textbf{3} (1972), 31--50. Univ. California Press,
Berkeley. MR0410936


\bibitem{JP74} Jain, N. C.; Pruitt, W. E.;
{\em Further limit theorems for the range of random walk}, J. Anal.
Math. \textbf{27} (1974), 94--117. MR0478361

\bibitem{Karamata}  Karamata, J.:
{\em Sur un mode de croissance r\'{e}guli\`{e}re},
Th\'{e}or\`{e}mes fondamentaux. (French) Bull. Soc. Math. France \textbf{61}
(1933), 55¨C62. MR1504998

\bibitem{LeGall86} Le Gall, J.-F.;
{\em Propri\'{e}t\'{e}s d'intersection des marches al\'{e}atoires.
I. Convergence vers le temps local d'intersection}, Comm. Math.
Phys. \textbf{104} (1986), 471--507. MR0840748

\bibitem{Revesz} R\'{e}v\'{e}sz, P.:
{\it Random Walk in Random and Non-Random Environments},
Second edition. World Scientific Publishing Co. Pte. Ltd., Hackensack, NJ, 2005.
xvi+380 pp. ISBN: 981-256-361-X. MR2168855

\bibitem{SV05} Salminen, Paavo; Vallois, Pierre;
{\em On First Range Times of Linear Diffusions}, J. Theoret. Probab.
\textbf{18} (2005), no. 3, pp. 567--593. MR2167642

\bibitem{Vallois96} Vallois, Pierre;
{\em The Range of a Simple Random Walk on $\Znum$}, Adv. in Appl.
Probab. \textbf{28} (1996), no. 4, pp. 1014--1033. MR1418244

\bibitem{VT97}  Vallois, Pierre; Tapiero, Charles S.
{\em Range Reliability in Random Walks}, Math. Methods Oper. Res.
\textbf{45} (1997), no. 3, pp. 325--345. MR1463263

\bibitem{WX12} Wu, Jun; Xie, Jian-Sheng;
{\em Range-Renewal Structure in Continued Fractions}, arXiv:1305.2088.

%\bibitem{XY} Xie, Jian-Sheng; Ying, Jiangang;
%{\em On the Range-Renewal of Random Sequence}, preprint.

 \end{thebibliography}



\end{document}